\documentclass[10pt]{amsart}

\usepackage{amsmath}
\usepackage{amscd}
\usepackage{amssymb}
\usepackage{latexsym}

\input xy
\xyoption{all}

\newtheorem{theorem}{Theorem}[subsection]
\newtheorem{lemma}[theorem]{Lemma}
\newtheorem{corollary}[theorem]{Corollary}
\newtheorem{proposition}[theorem]{Proposition}
\newtheorem{definition}[theorem]{Definition}
\newtheorem{axiom}[theorem]{Axiom}
\newtheorem{example}[theorem]{Example}
\theoremstyle{remark}
\newtheorem{remark}[theorem]{Remark}
\newcommand{\sn}{\par\smallskip\noindent}
\newcommand{\mn}{\par\medskip\noindent}
\newcommand{\bn}{\par\bigskip\noindent}

\newcommand{\cE}{\mathcal E}

\newcommand{\cM}{\mathcal M}

\newcommand{\cP}{\mathcal P}

\newcommand{\cT}{\mathcal T}
\newcommand{\cU}{\mathfrak U}

\newcommand{\hal}{{\rm H}}
\newcommand{\dS}{\mathcal S}
\renewcommand{\dim}{{\rm dim}\,}
\renewcommand{\ker}{{\rm ker}\,}

\newcommand{\Ass}{{\rm Ass}}
\newcommand{\ton}[1]{\stackrel{#1}{\to}}
\newcommand{\ev}{{\rm ev}}
\newcommand{\id}{1}
\newcommand{\rep}{{\rm rep}}
\DeclareMathOperator{\bimod}{-Bimod}
\DeclareMathOperator{\modi}{-Mod}
\DeclareMathOperator{\Hom}{HOM}
\DeclareMathOperator{\End}{End}
\DeclareMathOperator{\ind}{ind}
\newcommand{\Ab}{{\rm Ab}}
\newcommand{\map}{{\rm map}}
\newcommand{\colim}{\operatornamewithlimits{colim}}
\newcommand{\gaha}{{\Ass_\hal}}
\newcommand{\caha}{{\Ass^c_\hal}}
\newcommand{\aha}{{{\rm Alg}_\hal}}
\newcommand{\ac} {{{\rm Alg}_\hal^{\ind}}}
\newcommand{\acm}{{{\rm Alg}_\hal^\cM}}

\newcommand{\deq}{:=}
\newcommand{\sd}{\mathrm {sd}}
\newcommand{\sdi}{\sd^\bullet}

\newcommand{\ja}{\ell}
\def\Az{\mathbb{A}}
\def\Cz{\mathbb{C}}

\def\Nz{\mathbb{N}}

\def\Zz{\mathbb{Z}}
\def\Sz{\mathbb{S}}

\def\Spec{\rm{Spec}}
\def\Kt{\mathbf{K}}
\def\Kh{\mathbf{KH}}
\def\kaha{\Kh^{s}}
\def\Cc{\mathfrak{C}}
\def\Dc{\mathfrak{D}}
\def\Ac{\mathfrak{A}}

\newcommand{\weq}{\overset{\sim}{\to}}
\newcommand{\fib}{\twoheadrightarrow}
\title{Bivariant algebraic $K$-theory}
\author{Guillermo Corti\~nas}
\author{Andreas Thom}
\address{Guillermo Corti\~nas, Dep. Matem\'atica\\ Ciudad Universitaria Pab 1\\ (1428) Buenos Aires, Argentina\\ and  Dep. \'Algebra\\ Fac. de Ciencias\\
Prado de la Magdalena s/n\\ 47005 Valladolid, Spain.}
\email{gcorti@agt.uva.es} \urladdr{http://mate.dm.uba.ar/\~{}gcorti}
\thanks{Corti\~nas' research was partly supported by grants PICT03-12330, UBACyT-X294, VA091A05, and
MTM00958. Thom's research was partly supported by the EU-Network \textit{Quantum Spaces and Noncommutative Geometry} (Contract HPRN-CT-2002-00280)
and the DFG (SFB $478$ M\"unster, GK \textit{Metageometrie und analytische Topologie} M\"unster, GK \textit{Gruppen und
Geometrie} G\"ottingen).}
\address{Andreas Thom\\ Mathematisches Institut\\
Bunsenstr. 3-5\\ 37073 G\"ottingen, Germany.}

\email{thom@uni-math.gwdg.de}
\urladdr{http://www.uni-math.gwdg.de/thom/}

\begin{document}
\begin{abstract} We show how methods from K-theory of operator algebras can be applied in a completely algebraic setting to
define a bivariant, $M_\infty$-stable, homotopy-invariant, excisive K-theory of
algebras over a fixed unital ground ring $\hal$, $(A,B)\mapsto kk_*(A,B)$, which is universal in the sense that it maps uniquely to any other such theory.
It turns out $kk$ is related to C. Weibel's homotopy algebraic $K$-theory, $KH$. We prove that, if $\hal$ is commutative and $A$ is central as
an $\hal$-bimodule, then
\[
kk_*(\hal,A)=KH_*(A).
\]
We show further that some calculations from operator algebra $KK$-theory, such as the exact sequence of Pimsner-Voiculescu, carry over
to algebraic $kk$.
\end{abstract}
\maketitle

\tableofcontents

\section{Introduction}\label{sec:intro}
We consider associative, not necessarily unital or central algebras over a
fixed unital, not necessarily commutative ring $\hal$; we write $\gaha$ for the
category of such algebras. If $\hal$ is commutative, we consider also the category
$\caha$ of central algebras. Let $\aha$ be either of $\gaha$, $\caha$. Note that, by forgetting structure, we
can embed $\aha$ faithfully into each of the categories of
(central) $\hal$-bimodules, abelian groups and sets. Fix one of these
underlying categories, call it $\cU$, and let $F:\aha\to\cU$ be the
the forgetful functor. Let $\cE$ be the class of all exact
sequences of $\hal$-algebras
\[
(E)\quad 0\to A\to B\to C\to 0
\]
such that $F(B)\to F(C)$ is a split surjection. We construct a
triangulated category $kk$ with inverse suspension functor
$\Omega$, a functor $j:\aha\to kk$, and a family of connecting
maps
\[
\{\partial_E:\Omega(j(C))\to j(A)\}
\]
natural with respect to maps of exact sequences, such that the
following conditions are satisfied.
 \sn
\noindent{i)} For all $A\in
\aha$, $j(A\to M_\infty A)$ and $j(A\to A[t])$ are isomorphisms.
\sn
\noindent{ii)} For every exact sequence $E\in\cE$, the
following is a distinguished triangle in $kk$
\[
\xymatrix{\Omega
j(C)\ar[r]^{\partial_E}&j(A)\ar[r]&j(B)\ar[r]&j(C).}
\]
In particular, the sequence of functors
\begin{gather}\label{intro:kkn}
kk_n(A,B)\deq \hom_{kk}(j(A),\Omega^nj(B))\qquad (n\in\Zz)
\end{gather}
forms a bivariant, homotopy invariant, $M_\infty$-stable homology theory
which satisfies excision with respect to all sequences in $\cE$,
that is, for any $E\in \cE$ as above, $n\in\Zz$ and $D\in\aha$, we
have long exact sequences
\[
\xymatrix{\dots\ar[r]&kk_{n+1}(D,C)\ar[r]^{(\partial_E)_*}&kk_n(D,A)\ar[r]&kk_n(D,B)\ar[r]&kk_n(D,C)\ar[r]&\dots\\
          \dots\ar[r]&kk_{n+1}(A,D)\ar[r]^{(\partial_E)^*}&kk_n(C,D)\ar[r]&kk_n(B,D)\ar[r]&kk_n(A,D)\ar[r]&\dots}
\]
We call $kk_*$ {\it bivariant algebraic $K$-theory}. We show
further that $j$ is universal among all functors $X:\aha\to \cT$
with values in a triangulated category $\cT$, and equipped with
natural maps $\partial_E^X$ satifying the requirements above.
Universality means that if $X$ is such a functor, then there is a
unique triangulated functor $\bar{X}:kk\to \cT$, compatible with
connecting maps, such that $X=\bar{X}\circ j$. In particular
$\bar{X}$ induces, for each $n\in\Zz$, a group homomorphism
\begin{equation}\label{intro:bivaX}
kk_n(A,B)\to X_n(A,B)\deq \hom_{\cT}(X(A),\Omega^nX(B))\qquad
(n\in \Zz)
\end{equation}
compatible with all the structure. For example if $\hal$ is a
field of characteristic zero, and $\aha=\caha$, then applying this universal
property when $X$ is the Cuntz-Quillen pro-supercomplex \cite{cq},
we obtain a product preserving Chern character to bivariant
periodic cyclic cohomology
\[
ch_*:kk_*(A,B)\to HP^*(A,B).
\]
Bivariant periodic $K$-theory is related to C. Weibel's homotopy
algebraic $K$-theory, $KH_*$ (\cite{weih}). We prove that, for
commutative $\hal$, $\aha=\caha$, and $B\in\aha$, we have
\begin{equation}\label{kk=kh}
kk_*(\hal,B)=KH_*(B).
\end{equation}
It follows from this and the universal property of $kk$, that if
$X$ is any $M_\infty$-stable, homotopy invariant, excisive theory
as above, then there is a map
\[
KH_*(B)\to X_*(\hal,B).
\]
For example, when $\hal$ is a field of characteristic zero, we
obtain in this way the Chern character from $KH$ to periodic
cyclic homology $HP_*(B)=HP^*(\hal,B)$. \sn The construction of
$kk$ we present here is inspired in work of J. Cuntz on bivariant
$K$-theory of topological algebras; we generalize and adapt his
methods. In particular, no spectra are involved in the definition
of $kk$. On the other hand, one can also use the machinery of
symmetric spectra to construct a bivariant, homotopy invariant,
$M_\infty$-stable excisive theory as follows. Recall (e.g. from
\cite[Appendix D]{thom-thesis}) that $KH$ can be defined as the
homotopy groups of a certain functorial symmetric spectrum $\kaha:
\aha\to Sp^{\Sigma}$ with compatible external products
$\kaha(A)\land\kaha(B)\to\kaha(A\otimes B)$. Thus for each $A\in
\aha$, $\kaha(A)$ is a $\kaha(\hal)$-module spectrum. Let
$[,]_{\kaha(\hal)}$ denote the homomorphisms in the homotopy
category of $\kaha(\hal)$-module spectra. Put
\[
KH_n(A,B)\deq [\kaha(A),\Omega^n\kaha(B)]_{\kaha(\hal).}
\]
By the universal property of $kk$, there is a natural map
\begin{equation}\label{intro:compa}
kk_*(A,B)\to KH_*(A,B)
\end{equation}
This map is an isomorphism in some cases, for example when
$A=\hal$ and $\aha=\caha$; this yields \eqref{kk=kh}. Even for $\hal=\Zz$ (in which
case $\gaha=\caha$) the problem of determining for which $A$ and $B$ the map
\eqref{intro:compa} is an isomorphism is a difficult one.
This problem is discussed in Subsection \ref{subsec:boot}.
\mn
As said above, we adapt and generalize Cuntz' methods to the
algebraic setting; in particular continuous or differential
homotopy has to be replaced by algebraic (i.e. polynomial)
homotopy. In other words the interval $[0,1]$ has to be replaced
by the affine line $\Az^1=\Az^1_\Zz=\Spec(\Zz[t])$. This entails
some technical difficulties. For example, in the topological
setting, the basis for composing homotopies is that the interval
$I=[0,1]$ is homeomorphic to the amalgamated sum of two copies of
it; $I=[0,1/2]\cup [1/2,1]\cong I\vee I$. This is no longer true
in the algebraic setting;
$\Az^1\vee\Az^1=\Spec(\Zz[t]\times_\Zz\Zz[t])\not\cong \Az^1$. Then
there is also the fact that some of the homotopies used in the
operator and topological algebra setting to prove some key results are not polynomial,
so we need to come up with new algebraic homotopies to replace
them. All of these problems are dealt with in Section
\ref{sec:htpy}. Modulo these technicisms, the construction of $kk$
is carried out pretty much as in the topological algebra setting
(sections \ref{sec:ext}, \ref{sec:quasi}, \ref{sec:kkdef}), and
its universal property proved. Once the triangulated category $kk$
is constructed, some calculations performed by Cuntz in that
setting carry over with essentially the same proof
to the algebraic case. For example we show that if $A\to B$ and
$A\to C$ are split algebra monomorphisms, then $j$ maps the
coproduct $B\coprod_AC$ to the amalgamated sum $j(B)\oplus_{j
(A)}j(C)$ in $kk$
\begin{equation}\label{intro:pramalaga}
j(B\coprod_AC)\cong j(B)\oplus_{j(A)}j(C).
\end{equation}
In particular we get
\begin{equation}\label{intro:amalga}
kk_{n}(D,B\coprod_AC)=kk_n(D,B)\oplus_{kk_n(D,A)}kk_n(D,C) \qquad
(n\in \Zz),
\end{equation}
and similarly on the other variable. Moreover by the universal
property of $j:\aha\to kk$, it follows that if $X:\aha\to\cT$ is
any functor satisfying conditions i) and ii) above, then
\[
X(B\coprod_AC)\cong X(B)\oplus_{X(A)}X(C).
\]
In particular \eqref{intro:amalga} is also valid for the bivariant
groups $X_*$ of \eqref{intro:bivaX}.
As another example consider
an action $\Zz\to {\rm Aut}_\aha(A)$, sending $1\mapsto \alpha$,
and write $A\rtimes_\alpha \Zz$ for the crossed product algebra.
We show in \ref{subsec:cross} that there is a distinguished
triangle in $kk$
\begin{equation}\label{intro:pv}
\xymatrix{\Omega(j(A\rtimes_\alpha\Zz))\ar[r]&j(A)\ar[rr]^{1-j(\alpha^{-1})}&&j(A)\ar[r]&j(A\rtimes_\alpha\Zz)}
\end{equation}
Applying $kk_*(D,?)$ and $kk_*(?,D)$ we get a version in our
setting of the Pimsner-Voiculescu sequences of
\cite{PV} (see also \cite[14.3]{Cu_Weyl}). Again by the universal property, the same is
valid for any functor $X$ satisfying i) and ii) above.
Similarly, we prove that the Laurent polynomial ring is isomorphic
in $kk$ to the sum of $A$ and its suspension
\begin{equation}\label{intro:funda}
j(A[t^\pm])\cong j(A)\oplus \Sigma j(A).
\end{equation}
For example if $\hal$ is commutative and $\aha=\caha$, then applying $kk_*(\hal, )$ and using \eqref{kk=kh}, we
obtain the fundamental theorem of homotopy $K$-theory (\cite{weih})
\[
KH_n(A[t^\pm])=KH_n(A)\oplus KH_{n-1}(A)\qquad (n\in\Zz).
\]
Thus in light of the universal property of $kk$, we can interpret \eqref{intro:funda} as saying that any
$M_\infty$-stable, homotopy invariant, excisive theory satisfies the fundamental theorem.
\bn
The rest of this paper is organized as follows. In Section \ref{sec:conve} we fix some notations
used throughout the paper, and recall some basic facts about $\ind$-objects. Section \ref{sec:htpy} deals with polynomial homotopy. In \ref{subsec:alghtpy} we define the notion of homotopy between
($\ind$)-algebra homomorphisms and introduce an enrichment
of $\aha$ over the category $\Sz$ of simplicial sets. Thus for $A,B\in \aha$ we have a simplicial set $\hom_\aha^\bullet (A,B)$. We also introduce
a functor
\begin{equation}\label{intro:power}
\Sz^{op}\times \aha\to \aha,\ \ (X,A)\mapsto A^X
\end{equation}
such that
\[
\hom_{\Sz}(X,\hom_{\aha}^\bullet(A,B))=\hom_\aha(A,B^X)
\]
For example $B^{\Delta^n}=B\otimes \Zz[t_0,\dots,t_n]/\langle \sum t_i-1\rangle$. There is also a pointed version $B^{(X,\star)}$. The functor
\begin{equation}\label{intro:omega}
\Omega B\deq B^{(S^1,\star)}
\end{equation}
provides the loopspace (i.e. the inverse suspension) in the triangulated category $kk$ (see \ref{subsec:kktria}).
In \ref{subsec:subdi} we use \eqref{intro:power} to translate the subdivision of simplicial sets to algebras, and obtain a fibrant version
$\Hom_\aha^\bullet (A,B)$ of the enrichment above (see Theorem \ref{exinfi}). In \ref{subsec:pathomo} we define a path (or homotopy) between homomorphisms
$f_i:A\to B$, $i=0,1$ as a map $H:A\to B^{\sd^n\Delta^1}$ for some $n\ge 0$, such that for the evaluations we have $\ev_i H=f_i$. Equivalently,
a homotopy is an $\ind$-homomorphism from $A$ to the $\ind$-algebra
\[
B^{\sd^\bullet \Delta^1}\deq \{B^{\sd^n\Delta^1}\}_{n\ge 0}
\]
Similarly, a loop is an $\ind$-homomorphism to $B^{\dS^1}\deq B^{(\sd^\bullet S^1,\star)}$. We show in \ref{theo:pi1} that the set
of homotopy classes of loops can be recovered as the fundamental group of the fibrant simplicial set $\Hom^\bullet_\aha(A,B)$; we have
\[
[A,B^{\dS^1}]=\pi_1\Hom_{\aha}^\bullet(A,B).
\]
In Subsection \ref{subsec:rotatopy} we observe that a rotational
homotopy which appears in several proofs of topological algebra
$K$-theory can be implemented in the polynomial homotopy setting.
In Section \ref{sec:ext} we introduce a number of ($\ind$-)
algebras and $F$-split extensions which appear in the definition
of $kk$-theory and the $kk$-category. In \ref{subsec:matrix} we
introduce an $\ind$-algebra which will play the role the compact
(or rapid decay) operators play in topological setting. The
following subsections define the notion of extension  and
introduce algebraic versions of some extensions, such as the path,
universal, loop, and mapping path extensions. As in the
topological algebra setting, the universal extension classifies
all extensions up to homotopy (see \ref{class}). We also study (in
\ref{subsec:conext}) M. Karoubi's cone extension (\cite{KV1})
\begin{equation}\label{intro:conext}
M_\infty A\to \Gamma A\to \Sigma A.
\end{equation}
Interest in this extension comes from its use in delooping algebraic $K$-theory (\cite{KV1},\cite{Wagoner}). We show in Subsection \ref{subsec:conext} that
\eqref{intro:conext} is split as a sequence of $A$-bimodules, whence $F$-split. In \ref{subsec:infsum} we prove that $\Gamma A$
is an infinite sum ring in the sense of Wagoner \cite{Wagoner};
this is used further below, in Section \ref{sec:kkdef}, to prove that $\Gamma A$ is equivalent to zero in $kk$ and thus that the functor
$\Sigma$ is inverse to $\Omega$ in $kk$. Hence \eqref{intro:conext} is analogue to the operator algebra Calkin extension. The analogue of the operator
algebra Toeplitz extension is considered in \ref{subsec:toeplitz}. Section \ref{sec:quasi} is devoted to split exact and $M_2$-stable
functors defined on $\aha$ with values in an abelian category. We show in \ref{subsec:quasidefi} that $M_2$-stable functors are invariant under
(generalized) inner automorphisms (see \ref{prop:multiplier}), and observe that this implies any such functor is Morita invariant in the usual unital
sense (see \ref{morita}). In particular this will apply to $kk$-theory. Then we recall the notion of quasi-homomorphism and the extended functoriality properties of split-exact $M_2$-invariant functors
(see \ref{subsec:quasihomo}). Section \ref{sec:kkdef} is devoted to the definition and basic properties of $kk$. The section starts with the
definition of $kk(A,B)$; this makes use of the homotopy machinery developed in Section \ref{sec:htpy}. The composition product
$kk(A,B)\otimes kk(B,C)\to kk(A,C)$ and thus the category $kk$ are introduced in \ref{subsec:compopro}.
In \ref{subsec:exci} we prove excision in both variables of $kk$ (see \ref{theo:exci_1variable} and \ref{sequence}) and show that \eqref{intro:omega}
induces an equivalence in $kk$ (see \ref{isom-stab} and \ref{omegaisfunctor}). In \ref{subsec:deloop} we prove that the functor $\Sigma$ of
\eqref{intro:conext} is an equivalence in $kk$, which is pseudo-inverse to $\Omega$.
In \ref{subsec:kktria} we show that $kk$ is triangulated.
The universal property of $kk$ discussed above is
proved in \ref{subsec:uniprop}, where we also prove a second universal property of $kk$ (see \ref{univ:hom}). The latter says that any half exact, $M_\infty$-stable,
homotopy invariant functor $G$ from $\aha$ to an abelian category $\Ac$ factors uniquely through a homological functor
$\bar{G}:kk\to \Ac$. In \ref{subsec:compakas} we compare $kk$ with Kassel's bivariant $K$ group (see \cite{kas}). The latter is a defined for pairs
$(A,B)$ of unital algebras over commutative ground ring $\hal$. By definition, $K(\hal,B)=K_0(B)$ is the usual $K_0$; in particular,
$K(?,?)$ is not homotopy invariant. However, it is equipped with a composition product and a Chern character $K(A,B)\to HC^0(A,B)$
with values in (nonperiodic) cyclic homology. We show that there is a map $K(A,B)\to kk(A,B)$
compatible with all the structure.
Some computations, including \eqref{intro:pramalaga}, \eqref{intro:pv} and \eqref{intro:funda}, are carried out in  Section \ref{sec:compu}. We show
further (in \ref{subsec:gradednil}) that positively graded as well as nilpotent rings (satisfying a suitable $F$-splitting condition) go to zero in $kk$.
In Section
\ref{sec:compa} we compare $kk$ and $KH$. In \ref{subsec:khdef} we recall Weibel's definition of $KH$; we also give an alternative definition
which we prove is equivalent to the original one (see \ref{kh=kh'}), and from which a new proof of excision and nilinvariance of $KH$ is easily deduced
(see \ref{easy}).
This alternative definition is used in \ref{proof_main} to prove our main theorem \eqref{kk=kh}. The problem of determining when
\eqref{intro:compa} is an isomorphism is discussed in \ref{subsec:boot}.

\mn

\textit{Acknowledgements.} Many of the ideas that we are going to present go back to
the work of J.Cuntz. He analyzed Kasparov's $KK$-theory in
algebraic terms; this pointed to new ideas which finally could be
applied in many different settings. Our paper owes much to Cuntz'
pioneering work.
\sn
Part of the research for this article was carried out during
visits of the first author to the universities of G\"ottingen and M\"unster.
He is thankful to these institutions for their hospitality.
\mn
\section{Conventions and preliminaries}\label{sec:conve}
Let $\hal$ be a unital ring. The tensor product $\otimes_\hal$ makes the
category $\hal \bimod=\hal \otimes_\Zz \hal^{op}\modi$ of \hal-bimodules into
a monoidal category. By an {\it \hal-algebra} we understand a
monoid in $(\hal\bimod,\otimes_\hal)$, possibly without a neutral
element; thus algebras will be nonunital in general. We denote by
$\gaha$ the category of \hal-algebras. If $\hal$ happens to be commutative, we consider also the
full subcategory $\caha\subset\gaha$ of central algebras. Let $\aha$ be either of $\gaha$, $\caha$.
Most of our results do not depend on which of these two choices we pick, so they will be formulated
in terms of $\aha$. Whenever we need to restrict to one of the choices, this will be made clear.
In what follows an $\hal$-algebra will be an object of $\aha$.
\mn
Throughout, we assume fixed an underlying category $\cU$, which can be either the category of sets, of (central) bimodules
or of abelian groups. In each case we have a faithful forgetful functor $F: \aha \to \cU$ and a functor
$\tilde{T}: \cU \to \aha$, left adjoint to $F$.
\mn
If $\Cc$ is a category, we write $\ind-\Cc$ for the category of $\ind$-objects of $\Cc$.
It has as objects the directed diagrams in $\Cc$.
An object in $\ind- \Cc$ is described by a filtering partially ordered set
$(I,\leq)$ and a functor $X: I \to \Cc$. The set of homomorphisms from $(X,I)$ to $(Y,J)$ is
\begin{equation}
 \lim_{i \in I} \colim_{j \in J} \hom_{\Cc}(X_i,Y_j).
\end{equation}
There is a natural functor $\ind-(\ind-\Cc)\to \ind-\Cc$, mapping
\begin{equation}\label{collapse}
((X_i, J_i),I)\to (X_{ij},{\coprod}_{i\in I}J_i\times{i}).
\end{equation}
We shall use this functor to collapse any $\ind- \ind-$object to an
$\ind-$object. Any functor $F:\Cc\to \Dc$ extends to $\ind-\Cc\to
\ind-\Dc$ by $F(X_i,i\in I)=(F(X_i),i\in I)$. In particular any
functor $\Cc\to \ind-\Dc$ extends to $\ind-\Cc\to \ind-\ind-\Dc$,
which after collapsing gives a functor $\ind-\Cc\to \ind-\Dc$.

These extensions and collapsing shall be implicit, so that
whenever a functor is defined on a category we shall freely apply
it to the $\ind$-category. We shall identify objects of $\Cc$ with
constant $\ind$-objects, so that we shall view $\Cc$ as a
subcategory of $\ind-\Cc$.

\mn

Throughout, the letters $A,B,C, \dots$ denote ($\ind$-)$\hal$-algebras.
Thus $f: A \to  B$ is a homomorphism of ($\ind$-)$\hal$-algebras.

Let $L$ be a ring and $A$ be an $\hal$-algebra. For brevity,
we write $LA$ for the $\hal$-algebra $L \otimes_\Zz A$. Similarly, $AL$ stands for $A \otimes_\Zz L$.
\mn
We write $\Sz$ for the category of simplicial sets, and $\Nz$ for the set of natural numbers, which
for us include $0$. Thus $\Nz=\Zz_{\ge 0}$.
\mn
We use the notation $\deq$ to define the left side by the right side.

\section{Homotopy}\label{sec:htpy}

\subsection{Algebraic homotopies}\label{subsec:alghtpy}

Let $A$ be an \hal-algebra. Put \[ A^{\Delta^1} \deq A[t]=A \otimes_\Zz \Zz[t]. \]
There are several natural morphisms relating $A$ and $A[t]$. We
write $c_A: A \to A[t]$ for the inclusion of $A$ as constant
polynomials in $A[t]$ and $\ev_i: A[t] \to A$ for the evaluation
of $t$ at $i$ ($i \in \{0,1\}$). Note that $c_A$ is a section of
$\ev_i$.
\mn
Let $f_0,f_1: A \to B$ be morphisms in $\aha$. We
call $f_0$ and $f_1$ {\it elementary homotopic} if there exists a
morphism $h: A \to B[t]$ such that $\ev_i  h=f_i$ ($i\in\{0,1\}$).
Note that elementary homotopy is a reflexive and symmetric relation.
In general, it is not transitive.

\begin{definition}\label{def:htpy}
Let $f,g: A \to B$ be morphisms in $\aha$. We call $f$ and $g$ {\rm homotopic}, and write $f \sim g$, if they can be connected by a
chain of elementary homotopies. We denote the set of homotopy classes of morphisms from $A$ to $B$ by $[A,B]$. If
now $A=(A,I), \,B=(B,J)\in \ac$, we put
\[
[A,B]=\lim_j \colim_i[A_j,B_i].
\]
Note that there is a natural map $\hom_{\ac}(A,B)\to [A,B]$.
Two homomorphisms $f,g:A\to B$ in $\ac$ are called {\rm homotopic} if they have the same image in $[A,B]$.

\end{definition}

In order to organize the notion of algebraic homotopy, we introduce a simplicial enrichment of the category of \hal-algebras.
Note that
the assignment
$$
\Zz^{\Delta}: [n] \mapsto \Zz^{{\Delta}^n} \deq \Zz[t_0,\dots, t_n]/\langle 1-\sum_i t_i\rangle
$$
defines a simplicial unital ring. Let $A$ be an \hal-algebra; put
\[
A^{\Delta}:[n] \mapsto A\otimes_\Zz \Zz^{\Delta^n}.
\]
Using this construction, we can enrich the category $\aha$ over simplicial sets, as follows. We have
a mapping space functor $\hom^{\bullet}_{\aha}: (\aha)^{op} \times \aha \to \Sz$,
given by
\[
(A,B) \mapsto ([n] \mapsto \hom_{\aha}(A,B^{\Delta^n})).
\]
For $A,B,C\in \aha$, there is a simplicial map
\begin{equation}\label{compo}
\underline{\circ}:\hom^{\bullet}_{\aha}(B,C)\times \hom^{\bullet}_{\aha}(A,B)\to\hom^{\bullet}_{\aha}(A,C)
\end{equation}
which satisfies the axioms for simplicial composition \cite[I.1]{quimod}, so that $\aha$ equipped
with these data becomes a simplicial category in the sense of {\it loc.cit.} To define \eqref{compo}
we use the multiplication map $\mu: \Zz^{\Delta}\otimes\Zz^{\Delta}\to \Zz^{\Delta}$.
If $g\in \hom(B,C^{\Delta^n})$ and $f\in \hom(A,B^{\Delta^n})$, then
\[
g\underline{\circ}f:=(\id_C \otimes \mu) (g^{\Delta^n}\circ f)
\]
Here $g^{\Delta^n}$ is the map the functor $(?)^{\Delta^n}$ associates to $g$.
Furthermore, for every $A \in \aha$, the functor $\hom^{\bullet}_{\aha}(?,A): (\aha)^{op} \to \Sz$
has a left adjoint $A^?: \Sz \to (\aha)^{op}$. If $X\in\Sz$,
\begin{align*}
A^{ X} = &\lim_{\Delta^n \to X} A^{\Delta^n}\\
    = &\int_n\prod_{x\in X_n}A^{\Delta^n}.
\end{align*}
Here the first limit is taken over the category of simplices of $X$ (\cite[I.2]{Jardine}) and the integral sign
denotes an end \cite[Ch IX,\S5]{mac}. We have
\begin{equation}\label{adj1}
\hom_{\aha}(A,B^X)=\map_\Sz(X,\hom^\bullet_\aha(A,B)).
\end{equation}
Applying this to $A=t\hal[t]$ we obtain
\begin{equation}\label{map1}
B^X=\map_\Sz(X,B^\Delta).
\end{equation}
Let $\Sz_*$ be the category of pointed simplicial sets. For $(K,\star)\in\Sz_*$, put
\begin{align}\label{map2}
A^{{(K,\star)}}\deq&\map_{\Sz_*}((K,\star),A^\Delta)\nonumber\\
=&\ker(\map_\Sz(K,A^\Delta)\to \map_\Sz (\star,A^\Delta))\\
=&\ker(A^K\to A). \nonumber
\end{align}
\mn
\begin{lemma}\label{monepi}
Let $j:K\to L\in\Sz$, $\star\in K$, and $A\in\aha$. If $j$ is a cofibration, then $A^L\to A^K$ is surjective, and
the sequence
\[
 \xymatrix{0\ar[r]&A^{(L/K,\star)} \ar[r] & A^{(L,\star)} \ar[r] & A^{(K,\star)}\ar[r]&0 }
\]
is exact.
\end{lemma}
\begin{proof}
Under the identification \eqref{map1}, the map $A^L\to A^K$ is identified with
\[
j^*:\map_\Sz(L,A^\Delta)\to \map_\Sz(K,A^\Delta)
\]
Because $A^{\Delta}$ is a weakly contractible, fibrant simplicial set, $j^*$ is surjective
by the $LLP$ of cofibrations (\cite{quimod}). The remaining assertions of the lemma follow from the snake lemma
applied to the following map of exact sequences
\[
\xymatrix{0\ar[r]& A^{(L/K,\star)}\ar[d]\ar[r]&A^L\ar[d]\ar[r]&A^K\ar[d]\ar[r]&0\\
          0\ar[r]& 0 \ar[r]&A\ar[r]_1&A\ar[r]&0.}
\]
\end{proof}

\begin{proposition}\label{power=tensor}
Let $K$ be a finite simplicial set, $\star$ a vertex of $K$, and  $A$ an \hal-algebra. Then $\Zz^K$ and $\Zz^{(K,\star)}$ are
free abelian groups, and there are natural isomorphisms
\[
A\otimes\Zz^K\ton{\sim}A^K\qquad A\otimes\Zz^{(K,\star)}\ton{\sim}A^{(K,\star)}
\]
\end{proposition}
\begin{proof}
We consider the unpointed case first. The natural map
\[
\eta:A\otimes\Zz^K=A\otimes\lim_{\Delta^n\to K}\Zz^{\Delta^n}\to \lim_{\Delta^n\to K}A^{\Delta^n}
\]
is that induced by $A\otimes \Zz^{\Delta^n}=A^{\Delta^n}$. We shall prove by induction on $\dim K$ that if $K$ is finite
then $\Zz^K$ is free and $\eta$ is an isomorphism for all $A$. If $\dim K=0$ this is clear. Let $n\ge 0$ and assume both
assertions true for all finite simplicial sets of dimension $n$. If $K$ is finite and $\dim K=n+1$ we have a cocartesian square
\[
\xymatrix{{\coprod}_I\Delta^{n+1}\ar[r]&K\\
          {\coprod}_I \partial\Delta^{n+1}\ar[r]\ar[u]& sk^nK\ar[u]}
\]
where $I$ is a finite set. Applying the functor $A^?$ we get a cartesian square
\[
\xymatrix{{\prod}_IA^{\Delta^{n+1}}\ar[d]&A^K\ar[l]\ar[d]\\
          {\prod}_I A^{\partial\Delta^{n+1}}& A^{sk^nK}.\ar[l]}
\]
Both vertical arrows are surjective by \ref{monepi}. Moreover, because $I$ is finite, $\prod_I=\bigoplus_I$ in $\Ab$. Hence we have a short exact sequence of abelian groups
\begin{equation}\label{seqK}
0\to A^K\to A^{sk^nK}\oplus\bigoplus_IA^{\Delta^{n+1}}\to \bigoplus_IA^{\partial\Delta^{n+1}}\to 0.
\end{equation}
Applying this to $A=\Zz$, and taking into account that, by induction, $\bigoplus_I\Zz^{\partial\Delta^{n+1}}$ is free, we get a split exact
sequence. Thus $\Zz^K$ is a direct summand of a sum of one copy of $\Zz^{sk^nK}$, which is free by induction,
and of a finite number of copies of polynomial rings, which are also free. Hence $\Zz^K$ is free, and moreover, the sequence
\begin{equation}\label{seqtensK}
0\to A\otimes\Zz^K\to A\otimes\Zz^{sk^nK}\oplus A\otimes\bigoplus_I\Zz^{\Delta^{n+1}}\to A\otimes\bigoplus_I\Zz^{\partial\Delta^{n+1}}\to 0
\end{equation}

is exact. The natural map $\eta$ induces a map of exact sequences from \eqref{seqtensK} to \eqref{seqK}. By induction
this map is an isomorphism at both the middle and the right hand terms. It follows it is also an isomorphism at the
left. This proves the unpointed part of the proposition.
The part of the proposition concerning the pointed case follows from the unpointed one and the fact that the sequence of abelian groups
\[
0\to \Zz^{(K,\star)}\to \Zz^K\to \Zz\to 0
\]
is split.
\end{proof}

\begin{remark} The exponential law is not satisfied; in general
\begin{equation}
A^{ K\times  L} \not\cong (A^{ K})^{ L}
\end{equation}
Thus  $A^K$ is not an $A^K$-object in the sense of \cite[II.1, Def. 3]{quimod} (see  \cite[II.1, Prop. 1]{quimod}),
and therefore the axioms for a simplicial category in the sense of \cite[Def. 2.1]{Jardine} are not satisfied.
The failure of the exponential law already occurs when $K=\Delta^p$ and $L=\Delta^q$. Indeed,
\[
({A^{\Delta^p}})^{\Delta^q}=A^{\Delta^{p+q}}.
\]
On the other hand $\Delta^p\times \Delta^q$ is the amalgamated sum over $\Delta^{p+q-1}$ of $\binom{p+q}{q}$ copies of $\Delta^{p+q}$.
But since $A^?$ has a right adjoint, it maps colimits in $\Sz$ to colimits in $(\aha)^{op}$, that
is, to limits in $\aha$. In particular, $A^{\Delta^p\times \Delta^q}$ is the fiber product over
$A^{\Delta^{p+q-1}}$ of $\binom{p+q}{q}$ copies of $A^{\Delta^{p+q}}$.
For example
\[
A^{\Delta^1\times\Delta^1}=A^{\Delta^2{\coprod}_{\Delta^1}\Delta^2}=A^{\Delta^2}\times_{A^{\Delta^1}}A^{\Delta^2}\not\cong A^{\Delta^2}
\]
The reason for this is that $A^{\Delta^p}$ is really the ring of functions on the algebro-geometric affine space
$\mathbb{A}_\Zz^p$, and $\mathbb{A}_\Zz^p\times \mathbb{A}_\Zz^q=\mathbb{A}_\Zz^{p+q}$. Thus, with respect to products, affine spaces behave like
cubes, not simplices. This seems to indicate that an approach using cubical sets instead of simplicial sets would be appropriate,
but we do not follow this line.
\end{remark}

\subsection{Subdivision}\label{subsec:subdi}

Write $\sd: \Sz \to \Sz$ for the simplicial subdivision functor (see \cite[Ch. III.\S4]{Jardine}).
It comes with a natural transformation $h: \sd \to \id_{\Sz}$, which is usually called the last vertex map. We have an
inverse system
\[
\xymatrix{{\sd}^\bullet K:\sd^0K&&{\sd}^1K\ar[ll]_(.4){h_K}&&{\sd}^2K\ar[ll]_{h_{\sd K}}&&{\sd}^3K\ar[ll]_{h_{\sd^2 K}}&&{\dots}. \ar[ll]_
{h_{\sd^3 K}}}
\]
We may regard $\sdi K$ as a pro-simplicial set, that is, as an $\ind$-object in $\Sz^{op}$. The $\ind$-extension of the functor
$A^?:\Sz^{op}\to\aha$ maps $\sdi K$ to
\[
A^{\sdi K}=\{A^{\sd^n K}:n\in\Nz\}.
\]
If we fix $K$, we obtain a functor $(?)^{\sdi K}:\aha\to \ac$, which extends to
$(?)^{\sdi K}:\ac\to \ac$ in the usual manner explained in Section \ref{sec:conve}.

\mn
\begin{lemma}
Let $A\in\ac$. The functor $A^{\sdi ?}:\Sz^{op}\to \ac$ preserves finite limits.
\end{lemma}
\begin{proof}
Because $\sd$ has a right adjoint, it preserves all colimits. Similarly if $B$ is an algebra then $B^?$ maps colimits
in $\Sz$ to limits in $\aha$. Thus if $A=(A,I)\in\ac$, then
\[
A^{\sdi \colim_lK_l}=\{\lim_l A_i^{\sd^nK_l}:(i,n)\in I\times\Nz\}.
\]
To finish the proof it suffices to recall (e.g. from \cite[A.4]{am}) that finite limits in $\ac$ are computed levelwise.
\end{proof}

\mn

Let $A,B\in\ac$. The mapping space of Section \ref{subsec:alghtpy} extends to $\ind$-algebras by
\[
\hom^\bullet_{\ac}(A,B)\deq ([n]\mapsto \hom_{\ac}(A,B^{\Delta^n})).
\]
We also consider
\[
\Hom_{\ac}^\bullet(A,B)\deq ([n]\mapsto \hom_{\ac}(A,B^{\sdi\Delta^{n}})).
\]

\mn

\begin{proposition}\label{indadjoint}
Let $A$ and $B$ be $\ind-\hal$-algebras and let $K$ be a finite simplicial set. There is a natural isomorphism
\begin{equation*}
\map_\Sz(K,\Hom_{\ac}^\bullet(A,B)) = \hom_{\ac}(A, B^{\sdi K}).
\end{equation*}
\end{proposition}
\begin{proof}
The adjointness relation is checked as follows.
\[ \begin{array}{rcl}
\map_\Sz(K,\Hom_{\ac}^\bullet(A,B))
&=& \map_\Sz(\colim_{\Delta^n \to K} \Delta^n,\Hom_{\ac}^\bullet(A,B)) \\
&=& \lim_{\Delta^n \to K} \map_\Sz (\Delta^n, \Hom_{\ac}^\bullet(A,B)) \\
&=& \lim_{\Delta^n \to K} \hom_{\ac}(A,B^{\sdi\Delta^n})) \\
&=& \hom_{\ac}(A,\lim_{\Delta^n \to K} B^{\sdi\Delta^n})) \\
&=& \hom_{\ac}(A, B^{\sdi K}).\end{array}\]
\end{proof}

\begin{theorem}\label{exinfi}
Let $A\in \aha$, $(B,J)\in \ac$. Then
\[
\Hom_{\ac}^\bullet(A,B)=Ex^\infty\hom_{\ac}^\bullet(A,B).
\]
In particular, $\Hom_{\ac}^\bullet(A,B)$ is fibrant.
\end{theorem}
\begin{proof} The following chain of equalities is straightforward.
 \[\begin{array}{rcl} \map_\Sz(\Delta^k,\Hom_{\ac}^\bullet(A,B))
&=&  \colim_{(j,n) \in J \times \Nz} \hom_{\aha}(A,B_j^{\sd^n \Delta^k}) \\
&=&  \colim_{n\in \Nz} \colim_{j \in J} \map_{\Sz}(\sd^n \Delta^k, \hom_{\aha}^\bullet(A,B_j)) \\
&=&  \colim_{n\in \Nz} \map_{\Sz}(\sd^n \Delta^k, \colim_{j \in J}\hom_{\aha}^\bullet(A,B_j)) \\
&=&  \colim_{n\in \Nz} \map_{\Sz}(\Delta^k, Ex^n \colim_{j\in J} \hom_{\aha}^\bullet(A,B_j)) \\
&=& \map_\Sz (\Delta^k, Ex^{\infty} \hom_{\ac}^\bullet(A,B)).  \end{array} \]
\end{proof}

\subsection{Paths and homotopies}\label{subsec:pathomo}
Let $f_0,f_1:A\to B$ be \hal-algebra homomorphisms. A {\it path} from $f_0$ to $f_1$ is an $\ind$-homomorphism
$h:A\to B^{\sdi \Delta^1}$ such that $\ev_i(h)=f_i$ ($i=0,1$). Thus a path from $f$ to $g$ is
a $1$-simplex $\Hom_\ac^\bullet(A,B)$. Note that any path can be represented by a homomorphism $A\to B^{\sd^n\Delta^1}$ for some $n$, and conversely each such
map defines a path. We call a path from $f_0$ to $f_1$ {\it elementary} if it can be represented by a map $h:A\to B^{\Delta^1}$.
We use the notation $P(f_0,f_1)$ for the set of paths from $f_0$ to $f_1$. The
{\it constant path} of $P(f,f)$ is the degenerate path $s_0(f)$. We consider two a priori distinct notions of homotopy between paths. Two paths $h_0,h_1$ from $f$ to $g$ are called
{\it cubically homotopic} if there is an $\ind$-map $\kappa:A\to B^{\sdi {\Delta^1}}\otimes \Zz^{\sdi \Delta^1}$ such that
$(\ev_i\otimes 1)(\kappa)=s_0(f_i)$,
$(1\otimes \ev_i)(\kappa)=h_i$, ($i=0,1$). One checks that cubical homotopy of paths is an equivalence relation. We say that two paths $h_0,h_1$ as above are {\it simplicially
homotopic} if they are homotopic as $1$-simplices of $\Hom_\ac^\bullet(A,B)$ (\cite[3.1]{may}). We shall show that these two notions are equivalent.
Consider the free groupoid $\mathfrak{F}$ with set of objects $\hom_\aha(A,B)$ and one generating arrow $f\to g$ for each elementary
path from $f$ to $g$, where we identify a degenerate path $s_0(f)$ with the identity map of the object $f$. Thus
the following relation holds in $\mathfrak{F}$
\begin{itemize}
\item[i)] $s_0(f)=\id_f$ ($f\in\hom_{\aha}(A,B)$).
\sn
\noindent Write $\mathfrak{G}$ for the quotient of $\mathfrak{F}$ by the following relation
\sn
\item[ii)] If $\kappa:A\to B[t,u]=B^{\Delta^2}$ then $\kappa(?,1,t)\circ \kappa(?,t,0)=\kappa(?,t,1)\circ \kappa(?,0,t)$.
\sn
Consider also the groupoid $\mathfrak{G}'$ which results from modding out $\mathfrak{F}$ by the following relation

\item[ii')] If $\kappa:A\to B[t,u]$ then $\kappa(?,t,1-t)\circ \kappa(?,0,t) =\kappa(?,t,0)$.
\end{itemize}

\begin{lemma}
Let $A,B\in\aha$, $f,g\in\hom_\aha(A,B)$, and $h_1,h_2$ be paths from $f$ to $g$. Then $h_1$ and $h_2$ are cubically homotopic if and only if they
are simplicially homotopic. Moreover there is a groupoid structure with set of objects $\hom_\aha(A,B)$, where the maps $f\to g$ are
the homotopy classes of paths, and composition is given by concatenation. With this structure, $\hom_\aha(A,B)$ is isomorphic to both
$\mathfrak{G}$ and ${\mathfrak{G}'}$. In particular the latter two groupoids are isomorphic to each other.
\end{lemma}
\begin{proof}
That concatenation is compatible with homotopy is well-known in the simplicial case and straightforward in the cubical case.
Consider the map $\theta:P(f,g)\to \mathfrak{F}(f,g)$,
which sends $(h_0,\dots,h_{2n-1}):A\to B^{\sd^n\Delta^1}$ to the word $h_{2n-1}^{-1}\circ h_{2n-2} \dots\circ h_{1}^{-1}\circ h_0$. It
is clear that $\theta$ is surjective, and that it sends concatenation to composition.
Moreover we claim that $\theta(h)=\theta(h')$ in $\mathfrak{G}$ (respectively in $\mathfrak{G}'$) if and
only if $h$ and $h'$ are cubically (simplicially) homotopic. We prove the cubical part of the claim; a similar argument proves the simplicial
part. If $h$ and $h'$ are cubically homotopic paths then there is an $n\ge 1$ and a map $\kappa:A\to B^{\sd^n\Delta}\otimes\Zz^{\sd^n\Delta}$ such
that $(1\otimes\ev_0)\kappa$ represents $h$, $(1\otimes\ev_1)\kappa$ represents $h'$ and $(\ev_0\otimes 1)\kappa$ and $(\ev_1\otimes 1)\kappa$ represent respectively
$f\deq s_0(\ev_0(h))$ and $g\deq s_0(\ev_1(h'))$. If $n=1$ then, by relations i) and ii) above, $\kappa$ induces a commutative diagram in $\mathfrak{G}$
\[
\xymatrix{f\ar[r]^{\theta{h'}}&g\\
          f\ar[u]_1\ar[r]_{\theta{h}}&g\ar[u]^1}
\]
Hence $\theta(h)=\theta(h')$. For $n\ge 2$, $\kappa$ induces a larger square of arrows in $\mathfrak{G}$ made up of several commutative squares stuck together;
for example if $n=2$, we get
\[
\xymatrix{f\ar[r]^{\theta(h'_0)}&{\bullet}\ar[r]^{\theta(h'_1)^{-1}}&g\\
          f\ar[u]_1\ar[r] &{\bullet}\ar[r]\ar[u]&g\ar[u]^1\\
          f\ar[u]_1\ar[r]_{\theta(h_0)}&{\bullet}\ar[u]\ar[r]_{\theta(h_1)^{-1}}&g\ar[u]^1}
\]
It follows that the top arrow equals that at the bottom; this proves that $\theta$ maps cubically homotopic maps to equal maps in $\mathfrak{G}$.
Because $\theta$ maps concatenation to composition, to prove the converse it suffices to show that if $(h_0,\dots,h_{2n-1})$ represents a loop
$h$ based at $f$ and
\begin{equation}\label{rechloop}
h_{2n-1}^{-1}\circ h_{2n-2} \dots\circ h_{1}^{-1}\circ h_0=\id_f
\end{equation}
in $\mathfrak{G}$ then $h$ is cubically homotopic to $s_0(f)$. Now \eqref{rechloop} means that the left hand side is equal in $\mathfrak{F}$
to a composition of arrows of the form $h_1\circ h'_0\circ(h'_1\circ h_0)^{-1}$ where $h_i=(\ev_i\otimes 1)\kappa$ and
$h'_i=(1\otimes\ev_i)\kappa$ for some $\kappa:A\to B[t,u]$.
Note strings of concatenable elementary paths which differ by the intercalation of constant (i.e. degenerate) elementary paths represent
the same path. On the other hand any two strings without any configurations of the form
\begin{equation}\label{face2face}
\xymatrix{{{\bullet}}\ar[r]^h&{\bullet}&{\bullet}\ar[l]_h}
\end{equation}
which $\theta$ maps to the same arrow in $\mathfrak{F}$ represent the same path. Since $\kappa(?,t,u)=h(?,t)$ is a homotopy from \eqref{face2face}
to the constant path on $h(?,0)$, it suffices to show that pruning a string from any configurations
of the form
\[
\xymatrix{{\bullet}\ar[r]^{h'_0}&{\bullet} &{\bullet}\ar[l]_{\text{const.}}\ar[r]^{h_1}&{\bullet}&{\bullet}\ar[l]_{h'_1}\ar[r]^{\text{const}}&{\bullet}&{\bullet}\ar[l]_{h_0}}
\]
for $h_i$, $h'_i$ as above, does not change the path homotopy type. This is straightforward.
We show next that $\mathfrak{G}=\mathfrak{G}'$.
Let $\kappa$ be as in condition ii) above,
$h_i(?,t)=\kappa(?,i,t)$, $h'_i(?,t)=\kappa(?,t,i)$, $u(?,t)=\kappa(?,t,1-t)$. Put $\sigma:\Zz[t]\to\Zz[t]$, $\sigma(f)(t)=f(1-t)$.
Then for $\tilde{\kappa}\deq(1\otimes\sigma\otimes\sigma)\kappa$ we have $\tilde{\kappa}(?,0,t)=\sigma(h'_1)(?,t)$, $\tilde{\kappa}(?,t,0)=\sigma(h_1)(?,t)$ and
$\tilde{\kappa}(?,t,1-t)=\sigma(u)(?,t)$. By ii'), we have the following identities in $\mathfrak{G}'$:
\begin{equation}\label{rechnung}
u\circ h_0=h'_0\qquad \sigma(u)\circ \sigma(h_1)=\sigma(h'_1)
\end{equation}
On the other hand, if $h$ is any path, then $G(?,t,u)=h(?,u)$ satisfies $G(?,0,t)=h(?,t)$, $G(?,t,0)=h(?,0)$, $G(?,t,1-t)=\sigma(h)(?,t)$.
Hence $\sigma{h}$ is the inverse of $h$ in $\mathfrak{G}'$, by i) and ii). Thus it follows from \eqref{rechnung} that the relation
ii) holds in $\mathfrak{G}'$. Conversely, $\kappa'(?,t,u)\deq \kappa(?,t,(1-t)u)$ satisfies $\kappa'(?,0,t)=h_0$, $\kappa'(?,t,0)=h'_0$,
$\kappa'(?,1,t)=s_0(?,1,0)$ and $\kappa'(?,t,1)=u$. Hence relation ii') holds in $\mathfrak{G}$. We have proved that $\mathfrak{G}=\mathfrak{G}'$.
\end{proof}
Let $A\in\aha$. Define inductively
\[
 A^{\dS^1}\deq A^{(\sd^\bullet S^1,*)},\qquad A^{\dS^{n+1}}\deq (A^{\dS^n})^{\dS^1}.
\]
\begin{theorem}\label{theo:pi1}
Let $A\in\aha$, $B\in\ac$. There is a natural isomorphism
\[
[A, B^{\dS^1}] = \pi_1 \Hom_{\ac}^\bullet(A,B).
\]
\end{theorem}
\begin{proof} The left hand side is the group of automorphisms of the zero map in the groupoid $\hom(A,B)$ of the previous Lemma, whose arrows are
the cubical homotopy classes of paths. The right hand side is the set of loops based on the zero homomorphism modulo the simplicial homotopy
relation. By the lemma, both sides are isomorphic.
\end{proof}

\begin{remark} \label{eckmann} It is apparent that the different group structures on $[A, B^{\dS^n}]$ distribute over each other and share
a common neutral element. Therefore, by the Hilton-Eckmann argument, these group structures coincide and are abelian if $n \geq 2$.
\end{remark}

\subsection{Rotational homotopies}\label{subsec:rotatopy}
It was noted by several people in the field that most homotopies appearing in elementary proofs of
properties of $K$-theoretic invariants in operator algebras are either algebraic or rotational.
It is therefore not surprising that many ideas
of proofs, which were discovered in computations of the topological $K$-theory of operator algebras, have applications
to an algebraically homotopy invariant setting, once the rotational homotopies are under control.
\mn
For most of the purposes it is enough to construct an invertible matrix $W \in M_2(\Zz[t])$, such that
$$\ev_0(W) = \left( \begin{array}{cc} 1&0 \\0&1 \end{array} \right) \quad \mbox{and} \quad
\ev_1(W) = \left( \begin{array}{cc} 0&-1 \\1&0 \end{array} \right).$$

The matrix
$$W = \left( \begin{array}{cc} 1-t^2 &  t^3-2t \\ t & 1-t^2 \end{array} \right)$$
is a concrete example.

\begin{remark}
Note that this matrix $W$ is not orthogonal, in fact there is no orthogonal matrix with the
required properties.
Indeed, one easily checks that the existence of an orthogonal matrix with these properties
is equivalent to finding a solution to the equation
$A^2 + B^2 =1$ in $\Zz[t]$ such that $A_0=0, A_1=1, B_0=1$ and $B_1=0$. Each solution of $A^2+B^2=1$ in $\Zz[t]$ provides
an invertible element $A+Bi \in \Zz[i,t]$. However it is known that the only units in $\Zz[i,t]$ are $\pm 1$ and $\pm i$.
\end{remark}

\section{Extensions and Algebras}\label{sec:ext}

\subsection{Matrix algebras}\label{subsec:matrix}

Let $n$ be a positive integer; write $M_n$ for the algebra of $n \times n$ matrices with integer coefficients. If $n \leq m$, there is
a natural inclusion $\iota_{n,m}: M_n \to M_m$ of rings,
sending $M_n $ into the upper left corner of $M_m$. Thus the sequence $M_\bullet\deq(M_n)_n$ is an $\ind$-ring; write
$M_\infty$ for its colimit. We also consider the tensor product of $M_\bullet$ and $M_\infty$; this is the $\ind$-ring
$\cM_\infty=M_\bullet M_\infty$.

Let $A$ and $B$ be $\ind-\hal$-algebras. Put
\[
\{A,B\}\deq [A, \cM_\infty B].
\]
Here $[,]$ denotes homotopy classes of $\ind$-algebra homomorphisms, as defined in \ref{def:htpy}.
If $m\in\Zz^{\Nz\times\Nz}$ and $x\in M_\infty $ are matrices, then both $m\cdot x$ and $x\cdot m$ are well-defined elements of $\Zz^{\Nz\times\Nz}$.
The set
\begin{equation}\label{waggamma}
\Gamma^\ell\deq\{m\in\Zz^{\Nz\times\Nz}:m\cdot M_\infty\subset M_\infty\supset M_\infty\cdot m\}
\end{equation}
consists of those matrices in $\Zz^{\Nz\times\Nz}$ having
finitely many nonzero elements in each row and column. Consider the free abelian group
\[
\Zz^{(\Nz)}=\bigoplus_{n\in \Nz}\Zz.
\]
Note $\Gamma^\ell$ is a ring (called $\ell\Zz$ in \cite{Wagoner})
which is isomorphic to the subring of
$\End_\Ab(\Zz^{(\Nz)})$ consisting of those endomorphisms $f$ whose transpose maps $(\Zz^{*})^{(\Nz)}$ to itself. We remark
that if $T\in \Gamma^\ell$, then for every $p$ there exists an $n$ such that $T M_p\subset M_n\supset M_pT$. Thus
both left and right multiplication by $T$ define endomorphisms of the $\ind$-abelian group $M_\bullet$. In particular if $V,W\in\Gamma^\ell$
satisfy $VW=1$, then the rule $a\mapsto WaV$ defines ($\ind$-) ring endomorphisms $\psi^{V,W}$ of $M_\infty$ and
${\psi'}^{V,W}$ of $M_\bullet$.

\begin{lemma}\label{isomlema}
Let $V, W\in\Gamma^{\ell}$ be such that $VW= \id_{\Zz^{(\Nz)}}$, and $\psi=\psi^{V,W}$ and $\psi'={\psi'}^{V,W}$ as above.
Then $\iota\psi:M_\infty\to M_2M_\infty$
is homotopic to $\iota$, and $\psi':M_\bullet\to M_\bullet$ and $1\otimes\psi:\cM_\infty\to\cM_\infty$ are homotopic to the
identity maps.
\end{lemma}
\begin{proof} By \ref{subsec:rotatopy}, $\iota$ is homotopic to the inclusion $\iota'$ in the lower right corner. Thus
\[
\iota\psi^{V,W}=\psi^{V\oplus1,W\oplus 1}\iota\sim \psi^{V\oplus1,W\oplus 1}\iota' = \iota'\sim\iota.
\]
Hence $\iota(1\otimes\psi):\cM_\infty\to M_2\cM_\infty$ is homotopic to $\iota$. Moreover the same
calculation as above also shows that $\iota\psi'$ is homotopic to $\iota$. Hence it suffices to show that if $A$ and
$B$ are $\ind$-rings and $f,g:A\to M_\bullet B$ are homomorphisms such that $\iota f\sim\iota g$, then $f\sim g$. But this
is immediate from the defintion of the homotopy relation between $\ind$-homomorphisms and the fact that $\iota_{n,2n}$ is the
composite of $\iota\otimes 1_{M_n}$ and the isomorphism $M_2\otimes M_n\cong M_{2n}$.
\end{proof}

Next we recall two well-known operations in matrix rings. The first is the direct sum of matrices; this gives operations
$M_\infty\times M_\infty\to M_\infty$ and $M_\bullet\times M_\bullet\to M_\bullet$, each of which induces one on
$\cM_\infty$. It is straightforward from the previous lemma that these two agree up to homotopy, and are homotopy
associative and homotopy commutative, with the zero matrix as neutral element. Thus $\cM_\infty$ is a homotopy
abelian semigroup. The second operation is the tensor product of matrices. Again this gives two operations
$\star: M_\bullet\otimes M_\bullet\to M_\bullet$ and $\star: M_\infty\times M_\infty\to M_\infty$, which, together,
induce one on $\cM_\infty$, by the rule $(x\otimes y)\star (z\otimes w)=x\star z\otimes y\star z$. Again it is straightforward
from the previous lemma that, up to homotopy, the latter operation is associative, commutative, distributes over $\oplus$, and has
the $1\times 1$ matrix $e_{11}$ as neutral element. Hence $(\cM_\infty,\oplus,\star)$ is a homotopy semiring.
We shall consider the category $\acm$ whose objects are the $\ind-\hal$-algebras, and where the set of homomorphisms between two
objects $A,B$ is
\[
\hom_\acm(A,B) \deq \{A,B\}.
\]

\begin{theorem}\label{theo:homomatrix}
Let $A,B,C\in\ac$. There is an associative composition law
\[
\{B,C\} \times \{A,B\} \to \{A,C\}.
\]
which makes $\acm$ into a category, where the identity map of an object $A$ is the homotopy class of the inclusion
$A\to\cM_\infty A$. Moreover, direct sum of matrix blocks defines an enrichment of $\acm$ over abelian semigroups.
\end{theorem}
\begin{proof}
If $f:A\to \cM_\infty B$ and $g:B\to \cM_\infty C$ represent classes $[f]\in\{A,B\}$ and $[g]\in\{B,C\}$, define
the composite $[g]\star[f]$ as the class of the following composite map
\[
g\star f:A\stackrel{f}{\to} \cM_\infty B\stackrel{1\otimes g}{\longrightarrow}
\cM_\infty\cM_\infty C\stackrel{\star\otimes 1}{\longrightarrow}\cM_\infty C
\]
\end{proof}

Those $\ind-\hal$-algebras which are isomorphic in $\ac$ are called {\it matrix-homotopy equivalent}.
The assignment $[A,B] \ni f \mapsto f \otimes \iota_{\infty} \in \{A,B\}$ extends to a functor from the category
of $\ind-\hal$-algebras with homotopy classes of morphisms to the category $\acm$. Note that the induced map $[A,B] \to \{A,B\}$ is
an isomorphism if $B$ is \textit{stable}, i.e. isomorphic to $B \otimes \cM_\infty$.

\subsection{The notion of extension}\label{subsec:notion}

\begin{definition}
A sequence

\begin{equation}\label{extension}
\xymatrix{ A \ar[r]^f & B \ar[r]^g & C }
\end{equation}
of homomorphisms in $\ac$ is called an {\rm extension} if $f$ is a kernel of $g$ and $g$ a cokernel of $f$.
This implies that, up to isomorphism of arrows, $f$ and $g$ can be represented by level maps $f_\alpha$,$g_\alpha$ such that
$f_\alpha$ is kernel of $g_\alpha$ and $g_\alpha$ is a cokernel of $f_\alpha$ {\rm (\cite[A.3]{am})}.
\end{definition}

\begin{definition}
An extension of \hal -algebras \eqref{extension} is called {\rm $F$-split} if there exists a section
of $F(g)$ in $\cU$.
\end{definition}

\subsection{Path extension}\label{subsec:pathext}
Put
\[
\Omega=\Zz^{(S^1,\star)}
\]
Let $A$ be an \hal-algebra. By \ref{power=tensor} and \ref{monepi}, the
inclusion $\partial\Delta^1\subset \Delta^1$ induces an extension
\begin{equation}\label{basicsplit}
\xymatrix{\Omega A\ar[r]& A^{\Delta^1}\ar[r]^{(\ev_0,\ev_1)}&A\oplus A}
\end{equation}
Note this extension is naturally $F$-split.
An $\hal$-linear section of $(\ev_0,\ev_1)$ is $(a_0,a_1)\mapsto (1-t)a_0+ta_1$.

\subsection{Universal extension}\label{subsec:univext}

Consider the composite functor $T\deq\tilde{T} \circ F: \ac \to \ac$.
Note that, for $A\in\aha$, the counit map $\eta_A: T(A) \to A$ is
surjective (e.g. because $F$ is faithful \cite[IV.3, Thm 1]{mac}). Put
$J(A)\deq \ker\eta_A$. The following extension is called the {\it
universal extension} of $A$.
\[ \xymatrix{ J(A) \ar[r]^{\iota_A} &T(A) \ar[r]^{\eta_A} & A} \]

The name ``universal" is justified by the following observation.

\begin{proposition} \label{class}
Let $A \to B \to C$ be an $F$-split extension. There exists a commutative diagram of extensions as follows.
\[ \xymatrix{
A \ar[r] & B \ar[r] & C \\
J(C) \ar[u]^\xi \ar[r]^{\iota_C} & T(C) \ar[u] \ar[r]^{\eta_C} & C \ar[u]^{\id_C} } \]
Furthermore, $\xi$ is unique up to elementary homotopy. Because of this, we shall abuse notation and
refer to any such morphism $\xi$ as
{\rm the} classifying map of the extension.
\end{proposition}
\begin{proof} The existence of $\xi$ is clear from the adjointness of $F$ and $T$; its uniqueness up to homotopy
is straightforward from the fact that \eqref{basicsplit} is $F$-split.
\end{proof}
Note that if an extension is split in $\aha$, then its classifying map can be chosen to be zero.

\begin{proposition}\label{classmap}
Let \[ \xymatrix{
A  \ar[r] \ar[d]^f & B \ar[r] \ar[d] & C \ar[d]^g \\
A' \ar[r]        & B'\ar[r]        & C' } \]
be a commutative diagram of $F$-split extensions. Then there is a diagram
\[ \xymatrix{ J(C) \ar[d]^{J(g)} \ar[r] & A \ar[d]^f \\
J(C') \ar[r] & A'}\] of classifying maps, which is commutative up to elementary homotopy.
\end{proposition}

The following proposition clarifies the relation between classifying maps of extensions and tensor products.

\begin{proposition} \label{tensorcompatible}
Let $L$ be a ring, $A\in\aha$ and $\cU$ the underlying category. If $\cU$ is additive, then the extension
\begin{equation}\label{univL}
J(A)L\to T(A)L\to AL
\end{equation}
is $F$-split, and there is a choice for the classifying map $\phi_{A,L}: J(AL) \to J(A)L$ of this extension, which is
natural in both variables. If $\cU={\rm Sets}$ and $L$ is free as an abelian group, then \eqref{univL} is an extension,
and there is a choice of $\phi_{A,L}$ which is natural with respect to the first variable; in the second variable, it
is natural with respect to basis-preserving ring homomorphisms.
\end{proposition}
\begin{proof}
If $\cU$ is additive and $s:A\to T(A)$ the canonical splitting of the universal extension, then $s\otimes 1_L$
splits \eqref{univL} and is natural in both variables.
If $\cU={\rm Sets}$ and $L=\bigoplus_i\Zz$, then $\bigoplus_i s$ splits \eqref{univL}, and is natural in the first
variable, and also in the second if we restrict to basis-preserving homomorphisms. The proposition follows from this.
\end{proof}

The following corollary is easily implied by the preceding
proposition and Proposition \ref{power=tensor}.

\begin{corollary}\label{natural_map}
Let $K$ be a finite pointed simplicial set. There is a homotopy class of maps $J(A^K) \to J(A)^K$, natural with
respect to $K$,
which is represented by a classifying map of following the extension
\[
\xymatrix{ J(A)^K \ar[r]^{\iota_A^K} &T(A)^K \ar[r]^{\eta_A^K} & A^K.}
\]

Naturality means that for any  map of finite pointed simplicial sets $f: K \to L$, the following diagram commutes
up to homotopy

\[ \xymatrix{ J(A^K) \ar[r] \ar[d]^{J(A^f)} & J(A)^K \ar[d]^{J(A)^f} \\
J(A^L) \ar[r] & J(A)^L. }\]
\end{corollary}

In particular, it follows from the corollary that if $f,g: A \to A'$ and $f \sim g$, then $J(f) \sim J(g)$.
Moreover, Proposition \ref{classmap} and Corollary \ref{tensorcompatible} allow us to define the functor $J$ on the level of $\acm$. Indeed, given a homomorphism $f: A \to B \otimes \cM_\infty$
representing a homotopy class in $\{A,B\}$, we define $J(f) \in \{J(A),J(B)\}$ as the homotopy class of the composite
\[
\xymatrix{J(A) \ar[r]^(.4){J(f)}& J(B \otimes \cM_\infty) \ar[r]^{\phi_{B,\cM_\infty}}& J(B) \otimes \cM_\infty.}
\]

\subsection{Loop extension}\label{subsec:loop}
Write
\[
P \deq \Zz^{(\Delta^1,*)} \quad \mbox{and} \quad \cP \deq \Zz^{(\sd^\bullet \Delta^1,\star)}.
\]
If $A$ is an \hal-algebra, then by \ref{power=tensor} and \ref{monepi}, we have
a diagram of extensions
\[
\xymatrix{{\Omega} A\ar[d]_1\ar[r]&PA\ar[d]\ar[r]^{\ev_1}&A\ar[d]^{(0,1)}\\
          {\Omega} A\ar[r]& A^{\Delta^1}\ar[d]_{\ev_0}\ar[r]^{(\ev_0,\ev_1)}&A\oplus A\ar[d]^{\pi_1}\\
                             &A\ar[r]_1&A}
\]
Note that the rightmost column is split in $\aha$. Because \eqref{basicsplit} is naturally $F$-split,
both rows and the middle column are. In
particular
\begin{equation}\label{loopex}
\xymatrix{ {\Omega} A  \ar[r] & P A \ar[r]^{\ev_1} & A }
\end{equation}
is naturally $F$-split; we call it the {\it loop extension} of $A$.
Thus we can pick a natural choice for the classifying map, which we call
\[
\rho_A: J(A) \to \Omega A.
\]
A similar reasoning applies to the extension
\[ \xymatrix{ A^{\dS^1} \ar[r]& {\cP} A \ar[r]^{\ev_1} & A } \]

Since this is canonically split, it has a natural classifying map
$J(A) \to A^{\dS^1}$, namely the composite $J(A) \ton{{\rho}_A} \Omega A \ton{\iota} A^{\dS^1}$. We abuse
notation and call this map $\rho_A$ as well.

\subsection{Mapping path extension}\label{subsec:mapipath}
Let $f: A \to B$ be a morphism of \hal-algebras. The path extension of $f$ is the extension obtained from the loop extension
of $B$ by pulling it back to $A$:

\[ \xymatrix{
{\Omega} B \ar[r] \ar[d]^{\id_B}& P B \oplus_B A \ar[r] \ar[d] & A \ar[d]^f \\
{\Omega} B \ar[r] & P B \ar[r]^{\ev_1} & B. } \]

It has a natural $\hal$-linear section because the loop extension does. Its natural classifying map is %(up to strict homotopy)
$\rho_B \circ J(f): J(A) \to \Omega B$. We write $P_f\deq P B \oplus_B A$; $P_f$ is the {\it path-algebra} of $f$.
It comes with a natural evaluation map $\pi_f: P_f \to A$ and a natural inclusion
$\iota_f: \Omega B \to P_f$.
\mn
Again, similar reasoning applies in the subdivided setting, i.e. there is an extension
\[
\xymatrix{ B^{\dS^1} \ar[r]& {\cP} B \oplus_B A \ar[r]^(.6){\pi_f} & A }.
\]
We use the following notation
\[
\cP_f \deq \cP B \oplus_B A.
\]

\subsection{Cone extension}\label{subsec:conext}
Let $\Gamma$ be the ring of $\Nz \times \Nz$-matrices with
integer coefficients which satisfy the following two properties.
\begin{itemize}
\item[(i)] The set $\{a_{ij}, i,j \in \Nz \}$ is finite.
\item[(ii)] There exists a natural number $N \in \Nz$ such that each row and each column has at most
$N$ nonzero entries.
\end{itemize}

It is easily checked that $\Gamma$ is indeed a ring with respect
to entry-wise addition and the ordinary matrix multiplication law;
hence it is a subring of Wagoner's cone ring $\Gamma^\ell$ considered in \eqref{waggamma}.
In particular $M_\infty\subset \Gamma$ is an ideal; we put
\[
\Sigma=\Gamma/M_\infty.
\]
We note that $\Gamma$, $\Sigma$ are the cone and suspension rings
of $\Zz$ considered by Karoubi and Villamayor in
\cite[p.269]{KV1}, where a different but equivalent definition is
given. In keeping with our notational conventions, we shall write
$\Gamma A$ and $\Sigma A$ for $\Gamma\otimes A$ and $\Sigma
\otimes A$. As in \cite{KV1}, one may also consider the ring
$\Gamma(A)$ of matrices with entries in $A$ satisfying i) and ii)
above, and the quotient $\Sigma(A)=\Gamma(A)/M_\infty A$. There
are natural ring homomorphisms
\begin{equation}\label{gammagamma}
\Gamma A\to \Gamma(A),\qquad \Sigma A\to \Sigma(A)
\end{equation}
\begin{lemma} The maps \eqref{gammagamma} are isomorphisms.
\end{lemma}
\begin{proof} It suffices to show that $\Gamma A\to \Gamma(A)$ is an isomorphism. Call this map $f$.
For the purpose of this proof, $\Nz\times\Nz$-matrices will be regarded as maps
on $\Nz\times\Nz$. Write $\pi_i:\Nz\times\Nz\to\Nz$ for each of the projection maps.
The characteristic function $\chi_S$ of a subset $S\subset \Nz\times\Nz$ is in $\Gamma$ if and only
if
\[
S\in\Theta\deq\{T\subset\Nz\times\Nz: (\exists N)(\forall n\in\Nz), \#(\pi^{-1}_i(n)\cap T)\le N, \ \ i=1,2\}
\]
If $a\in \Gamma A$, then
\begin{align*}
a=&\sum_{\lambda\in {\rm Im} (a)}\lambda\cdot \chi_{\{(i,j):a_{ij}=\lambda\}}\\
=&f(\sum_{\lambda\in {\rm Im} (a)}\lambda \otimes\chi_{\{(i,j):a_{ij}=\lambda\}})
\end{align*}
Hence $f$ it surjective. In particular $\Gamma$ is generated as an abelian group by the characteristic functions of the elements of
$\Theta$. Thus every element $x$ of $\Gamma A$ can be written as a sum
\begin{equation}\label{equis}
x=\sum_{i=1}^na_i\otimes\chi_{S_i}
\end{equation}
for some $S_1,\dots,S_n\in\Theta$ with $S_i\neq S_j$ if $i\neq j$. Put
\[
S=\bigcup_{i=1}^nS_i.
\]
For each $F\subset\{1,\dots,n\}$, write
\[
S_F=\bigcap_{i\in F} S_i\cap\bigcap_{j\notin F} (S\backslash S_j)
\]
Note the $S_F$ are disjoint from one another, and are elements of $\Theta$. Hence
\[
\chi_{S_i}=\sum_{i\in F}\chi_{S_F}.
\]
Therefore the element \eqref{equis} can be rewritten as
\[
x=\sum_F(\sum_{i\in F}a_i)\otimes\chi_{S_F}.
\]
 Thus if $f(x)=0$, then $\sum_{i\in F}a_i=0$ for all $F$, whence $x=0$.
\end{proof}

An important property of the cone ring $\Gamma$ is that the sequence
\begin{equation}\label{noebe}
0\to M_\infty\to\Gamma\to \Sigma\to 0
\end{equation}
is a split exact sequence of free abelian groups. This follows from the theory of Specker groups developed
by  G. N\"obeling in \cite{noeb}. Almost by definition $M_\infty \subset \Gamma$ is an inclusion of Specker
subgroups of the abelian group of bounded maps $\Nz \times \Nz \to \Zz$ (see \cite{noeb}).
By \cite[Theorem 2]{noeb}, this implies that \eqref{noebe} is a split-exact sequence of free abelian groups.
In particular, it follows that for every $\hal$-algebra $A$,
\begin{equation}\label{conext}
\xymatrix{ M_\infty A  \ar[r] & \Gamma A \ar[r] & \Sigma A }
\end{equation}
is an extension, and is split as a sequence of $\hal$-bimodules. \mn

The rings $\Gamma$ and $\Gamma^\ell$ share several formal properties. For example, Wagoner showed that $\Gamma^\ell$ is an infinite sum
ring; the same is true of $\Gamma$, as we shall see in Lemma
\ref{infinitesum}.

\begin{remark}
We remark that there is a formally similar definition of $\Gamma$
in \cite[Sec. 10]{Wodz} which in fact does {\bf not} define a ring. The definition in \cite[Sec. 10]{Wodz} relaxes condition
(ii)
by requiring only that the number of nonzero elements in each row and column be finite.
The resulting set of matrices is not closed under products. To see this,
consider the matrix $A$ given by direct sum of matrix blocks $A_n$, ($n \in \Nz_{\ge 1}$), such that each block
$A_n\in M_n \Zz$ has all its entries equal to $1$. The matrix $A$ satisfies (i) and has finitely many nonzero entries
in each row and column; however $A^2$ does not satisfy (i).
\end{remark}

\subsection{Infinite sum rings}\label{subsec:infsum}

The notion of sum rings and infinite sum rings was introduced by Wagoner in \cite{Wagoner}. In the context of operator algebras,
similar algebras were introduced and sucessfully studied by Cuntz \cite{Cuntz2}.

\begin{definition}
\sn
\begin{itemize}
\item
A {\rm sum $\hal$-algebra} is an algebra $A\in\aha$ together with elements $\alpha_1,\alpha_2, \beta_1,\beta_2$,
which satisfy the relations
\begin{gather*}
\alpha_1\beta_1 = \alpha_2\beta_2 =1, \quad \beta_1 \alpha_1 + \beta_2 \alpha_2 = 1\\
[\alpha_i,h] = [\beta_i,h] =0, \quad (h\in \hal,\ \ i=1,2).
\end{gather*}
\item
Let $a$ and $b$ be elements in a sum \hal-algebra. We use the notation \[ a \oplus b = \beta_1 a \alpha_1 + \beta_2 b \alpha_2. \]
\item
Let $B$ be an \hal-algebra and let $\phi, \psi: B \to A$ be morphisms in $\aha$ into a sum \hal-algebra $A$.
We write $\phi \oplus \psi$ for the \hal-algebra homomorphism which sends \[B \ni b \mapsto \phi(b) \oplus \psi(b) \in A. \]
\item
An {\rm infinite sum $\hal$-algebra} is a sum $\hal$-algebra $A$ together with a $\hal$-algebra
endomorphism $\phi^{\infty}: A \to A$ which satisfies
\[
\id_A \oplus \phi^{\infty} = \phi^{\infty}.
\]
\end{itemize}
\end{definition}

\begin{lemma} \label{infinitesum} Let $A$ be a unital \hal-algebra. Then $\Gamma A$ is an infinite sum \hal-algebra. \end{lemma}
\begin{proof} First of all, $\Gamma A$ is a sum \hal-algebra, with
\[
\alpha_1 =\sum_{i=0}^\infty e_{i,2i}, \quad \beta_1 =\sum_{i=0}^\infty e_{2i,i},  \quad
\alpha_2 =\sum_{i=0}^\infty e_{i,2i+1}, \quad \mbox{and} \quad \beta_2 =\sum_{i=0}^\infty e_{2i+1,i}.
\]
Let $a\in\Gamma A$. Because the map $\Nz\times\Nz\to \Nz$, $(k,i)\mapsto 2^{k+1}i+2^k-1$, is injective, the following assignment
defines gives a well-defined, $\Nz\times\Nz$-matrix
\[
\phi^\infty(a) = \sum_{k=0}^\infty \beta_2^k \beta_1 a \alpha_1 \alpha_2^k=\sum_{k,i,j}e_{2^{k+1}i+2^k-1,2^{k+1}j+2^k-1}\otimes a_{i,j}.
\]
One checks that in fact $\phi^\infty(a)\in\Gamma A$. Next note that $\alpha_2 \beta_1= \alpha_1 \beta_2 =0$. It follows that
$\alpha_1 \alpha_2^i \beta_2^j \beta_1 = \delta_{i j}$ and hence that $\phi^{\infty}$ is a homomorphism of algebras.
The conditions in the definition of infinite sum algebra are straightforward.
\end{proof}
The following lemma is recalled from \cite{Wagoner}.
\begin{lemma}{\rm (Wagoner, \cite[page 355]{Wagoner})}\label{wagoner_lemma} Let $A$ be a sum \hal-algebra. There is an invertible
$3 \times 3$-matrix $Q \in M_3A$ which satisfies
\[ Q \left( \begin{array}{ccc} a \oplus b &0&0 \\ 0&0&0\\0&0&0 \end{array} \right) Q^{-1} =
\left( \begin{array}{ccc} a &0&0 \\ 0&b&0\\0&0&0 \end{array} \right), \] for
all $a,b \in A$.
\end{lemma}

\subsection{Laurent extension}\label{subsec:laurent}

Write $\Zz[t^{\pm 1}]$ for the Laurent polynomial ring. It comes with a natural evaluation map $ev_1: \Zz[t^{\pm 1}] \to \Zz$ which sends
$t$ to $1$. This surjection is split and defines a split extension as follows.

\[ \xymatrix{ {\sigma} A  \ar[r] & A[t^{\pm 1}] \ar[r]^{ev_1} & A }. \]
Here, $\sigma = (t-1)\Zz[t^{\pm 1}] = \ker(\ev_1)$.
\subsection{Toeplitz extensions}\label{subsec:toeplitz}

Let $\tau_\hal$ be the \hal-algebra which is unital and free on two generators $\alpha$ and $\beta$ satisfying the relations
$\alpha \beta=1$ and $[\alpha,h]= [\beta,h] =0$, for $h \in \hal$. It is easy to verify that $\tau_\hal \cong \tau_\Zz \otimes_\Zz \hal$, as
\hal-algebras. Where no confusion can arise, the subscript $\hal$ shall be omitted.

The natural map $\tau_{\Zz} \to \Zz[t^{\pm 1}]$, mapping $\alpha$ to $t$ and
$\beta$ to $t^{-1}$, is surjective and its kernel identifies with $M_{\infty}$ as the following lemma shows. It also follows from classical results
of N. Jacobson in \cite{Jacobson}.

\begin{lemma}
The kernel of the natural map $\tau \to \Zz[t^{\pm 1}]$ is isomorphic to $M_\infty$.
\end{lemma}
\begin{proof} The elements $\hat{\alpha}=\sum_ie_{i,i+1}$, $\hat{\beta}=\sum_ie_{i+1,i}$ of $\Gamma$ satisfy
$\hat{\alpha}\hat{\beta}=1$. Write $\hat{}$ for the ring homomorphism $\tau\to \Gamma$, $\alpha\mapsto \hat{\alpha}$,
$\beta\mapsto\hat{\beta}$. Since $1-\hat{\beta}\hat{\alpha}=e_{11}\in M_\infty$, $\hat{}$\ \ induces a map of extensions
\[
\xymatrix{{\ker} \pi\ar[r]\ar[d]&{\tau}\ar[d]\ar[r]^\pi&{\Zz}[t^{\pm}]\ar[d]\\
           M_\infty \ar[r]&{\Gamma}\ar[r]&\Sigma}
\]
As an abelian group, $\tau$ is generated by $\{\beta^p\alpha^q:p,q\ge 0\}$. Since
the elements $\hat{\beta}^p\hat{\alpha}^q=\sum_je_{j+p,j+q}$ are linearly independent, it follows
that $\tau\to\Gamma$ is injective. In particular $\ker\pi\to M_\infty$ is injective. On the other hand
\[
\hat{\beta}^p\hat{\alpha}^q-\hat{\beta}^{p+1}\hat{\alpha}^{q+1}=e_{p,q}
\]
This proves that $\ker\pi\to M_\infty$ is surjective, whence an isomorphism.
 \end{proof}

Let $A$ be an \hal-algebra. From the lemma above, we get an extension
\[ \xymatrix{ M_\infty A  \ar[r] & \tau A \ar[r] & A[t^{\pm 1}],} \] and an analogous extension
\[ \xymatrix{ M_\infty A  \ar[r] & \tau_0 A \ar[r] & {\sigma} A.} \]
In the last line, $\tau_0 = \tau \oplus_{A[t^{\pm 1}]} \sigma A$.

\subsection{The algebras $q(A)$ and $Q(A)$}\label{subsec:qAQA}

Let $A$ be an \hal-algebra; write $Q(A)$ for the coproduct of $A$ with itself in the category $\aha$. Note that there
is a natural codiagonal map $m: Q(A) \to A$. Write $q(A)\deq \ker m$. The algebra $q(A)$ was introduced by J.Cuntz in
his study of $KK$-theory, see e.g. \cite{Cuntz1}. We have the following extension:

\[ \xymatrix{ q(A)  \ar[r] & Q(A) \ar[r] & A.} \]
Note this extension is split in $\aha$ in two different ways, corresponding to each of the canonical inclusions $A\to Q(A)$. Moreover,
one checks that it is universal among such bi-split extensions.

\section{Split-exact and $M_2$-stable functors}\label{sec:quasi}

\subsection{Definitions}\label{subsec:quasidefi}

In this subsection we recall some basic facts about split-exact
functors. Most of the results are well-known in the setting of
topological algebras (see \cite{Cu_lecturenotes}) and adapt to the
algebraic setting in a straightforward way. \mn
A functor $E$ from
$\aha$ to an abelian category is called {\it split-exact}, if it sends split-extensions in
$\aha$
to split-extensions in the abelian category. A functor $E$ from
$\aha$
to an abelian category is called {\it$M_2$-stable}, if
for every $A\in \aha$, the map $E(A) \to E(M_2A)$ induced by the natural inclusion $\iota: A \to M_2A$ is an isomorphism.

\begin{proposition} \label{basic}
Let $E$ be a split-exact functor and $f,g: A \to B$ morphisms in $\aha$.
If $f(a)g(b)=g(a)f(b)=0$, for all
$a,b \in A$, then $(f+g): A \to B$ is a morphism in $\aha$ and $E(f) + E(g) = E(f+g)$.
\end{proposition}

The proof of the preceding proposition is an easy consequence of split-exactness.

\begin{proposition}\label{prop:multiplier} Let $E$ be an $M_2$-stable functor, $B$ an $\hal$-algebra, $A\subset B$ a subalgebra, and
$V,W\in B$ elements such that
\[
WA, VA\subset A,\ \ aVWa'=aa'\ \ (a,a'\in A),\ \ [V,H]=[W,H]=0.
\]
Then
\[
\phi^{V,W}:A\to A,\quad a \mapsto WaV
\]
is an $\hal$-algebra homomorphism, and
\[
E(\phi^{V,W}) = \id_{E(A)}.
\]
\end{proposition}
\begin{proof}
First of all, note that the inclusions
\[ a \mapsto \left( \begin{array}{cc} a & 0 \\ 0 & 0 \end{array} \right) \quad \mbox{and} \quad
a \mapsto \left( \begin{array}{cc} 0 & 0 \\ 0 & a \end{array} \right) \]
induce the {\it same} isomorphism $E(A) \to E(M_2 A)$. This is true, since on $M_4 A$, there are
automorphisms of order $2$ and $3$, mapping the different inclusions onto each other; i.e. $M_2$-stability is
used once more.

The argument is finished (as in the proof of Proposition \ref{isomlema})
by noting that $\phi^{V \oplus \id,W \oplus \id}$ restricts to $\phi^{V,W}$ in the upper left corner and to $\id_A$ in the lower right corner.
\end{proof}
\begin{remark}\label{morita}
Note that $M_2$-stable functors are $M_n$-invariant for all $n$. Further, we shall presently use \ref{prop:multiplier} to show that they are also Morita
invariant. Let $A,B\subset C$ be subalgebras of an $\hal$-algebra $C$, $n\ge 1$ and $p,q\in C^n$ such that $[p,H]=[q,H]=0$,
 $p_iA, Aq_i\subset B$ ($1\le i\le n$) and $a(\sum_iq_ip_i)a'=aa'$ ($a,a'\in A$). Then
\[
\xi_{p,q}:A\to M_nB,\ \ a\mapsto e_{i,j}\otimes p_iaq_j
\]
is an $\hal$-algebra homomorphism. Thus if $E$ is an $M_2$-stable functor, we have a homomorphism
\[
E(p,q)\deq E(\iota_n)^{-1}E(\xi_{p,q}):E(A)\to E(M_nB)\cong E(B).
\]
In the particular case when $A=B$, for $V=\sum_ie_{1,i}\otimes q_i$ and $W=\sum_i e_{i,q}\otimes p_i$ we have
$\xi_{p,q}=\phi^{V,W}\circ \iota_n$. Thus in this case $E(p,q)$ is the identity map, by \ref{prop:multiplier}.
Now we apply this to show that if $A$ and $B$ are Morita equivalent unital algebras, then $E(A)\cong E(B)$.
Suppose the equivalence is implemented by $\hal$-central bimodules
${}_BP_A$,${}_AQ_B$ and isomorphisms $f:P\otimes_A Q\to B$ and $g:Q\otimes_B P\to A$. Let
\[
C=\left[\begin{matrix}A & Q\\ P & B\end{matrix}\right].
\]
Choose $n$, $p,p'\in P^n$, $q,q'\in Q^n$ such that
$f(\sum_ip_i\otimes q_i)=1$, $g(\sum_iq'_i\otimes p'_i)=1$. Then
by what we have just seen, $E(p,q)$ and $E(q',p')$ are inverse
isomorphisms; in particular $E(A)\cong E(B)$. Hence $E$ is Morita
invariant. Note that this gives an obvious way of introducing a
notion of Morita invariance for nonunital $\hal$-algebras which
will be automatically preserved by $M_2$-invariant functors (see
\cite[\S7]{Cu_Weyl} for a similar notion in the topological
algebra context).
\end{remark}
\subsection{Quasi-homomorphisms}\label{subsec:quasihomo}

The notion of quasi-homomorphism is a formalism which helps organize the extended functoriality properties of split-exact functors.
We summarize some basic results from \cite{CuntzThom} and \cite{Cu_lecturenotes}.

\begin{definition}
Let $A$ and $B$ be \hal-algebras. A {\rm quasi-homomorphism} from $A$ to $B$ is a pair of \hal-algebra homomorphisms
$\phi, \psi: A \to D$ into an \hal-algebra $D$, which contains $B$ as an ideal, and is such that
\[
\phi(a)-\psi(a)\in B\qquad (a\in A).
\]
We use the notation
\[
(\phi,\psi) : A \to D \rhd B.
\]
\end{definition}

\begin{example}
The canonical inclusions $\iota_i:A\to Q(A)$, $i=1,2$ define a quasi-homomorphism $(\iota_1,\iota_2): A \to QA \rhd qA$.
This quasi-homomorphism is
universal in the sense that for any quasi-homomorphism $(\phi,\psi) : A \to D \rhd B$ there is a unique homomorphism
$\phi\ast \psi:QA\to D$ making the obvious diagram commute. The induced map
\[
\eta:q(A)\to B
\]
is called the {\rm classifying map} of $(\phi,\psi)$.
\end{example}
The following theorem sumarizes the extended functoriality of split-exact functors under quasi-homomorphisms. For a proof,
see \cite[3.2]{CuntzThom}.

\begin{theorem} \label{quasi-homomorphisms} Let $E$ be a split-exact functor from $\aha$
to an abelian category. Furthermore, let $(\phi,\psi): A \to D \rhd B$ be a quasi-homomorphism. Then, there is an induced map
\[ E(\phi,\psi): E(A) \to E(B) \] which satisfies the following naturality conditions:
\begin{itemize}
\item If $\psi=0$, then $E(\phi,0)=E(\phi)$.
\item $ E(\phi,\psi) = - E(\psi,\phi) $
\item If $\alpha: C \to A$ is an \hal-algebra homomorphism, then
\[ E(\phi  \alpha,\psi  \alpha) = E(\phi,\psi)  E(\alpha). \]
\item If $\eta: D \to D'$ is an \hal-algebra homomorphism which maps $B$ into an
ideal $B' \lhd D'$, then
\[ E(\eta  \phi,\eta  \psi) = E(\eta|_B)  E(\phi,\psi).\]
\end{itemize}
\end{theorem}

The following corollary clarifies in what sense $(\iota_1,\iota_2): A \to Q(A) \rhd q(A)$ is the universal quasi-homomorphism.

\begin{corollary} Let $(\phi,\psi): A \to D \rhd B$ be a quasi-homomorphism, $\eta:q(A)\to B$ its classifying map, and $E$
a split-exact functor. Then $E(\phi,\psi) = E(\eta) \circ E(\iota_1,\iota_2)$.
\end{corollary}

A proof of the following lemma can be found in \cite{CuntzThom}.

\begin{lemma} \label{split_sum}
Let $E$ be a split exact functor. If $\phi,\psi: A \to B$ satisfy $\phi(a)\psi(b)=\psi(a)\phi(b)=0$, for all $a,b \in A$,
then $E(\phi + \psi, \phi) = E(\psi)$.
\end{lemma}

\section{Algebraic $kk$-theory}\label{sec:kkdef}

\subsection{Definition of $kk$-theory}\label{subsec:kkdef}
Let $A$ and $B$ be $\ind-\hal$-algebras. Let $f: A \to B$ be a morphism in $\ac$. There is an associated map
$J(A) \to B^{\dS^1}$ which is given by $\rho_{B}  J(f)$. The map $\rho_{B}J(f)$ is a classifying map for the
mapping path extension
\[ \xymatrix{ B^{\dS^1} \ar[r] & {\cP_f} \ar[r]^{\pi_f}& A. } \]
\mn
Let us now assume that $A$ and $B$ are constant, i.e. just \hal-algebras. We set:
\[ kk(A,B) = \colim_{n \in \Nz}\{J^n(A),B^{\dS^n}\} = \colim_{n \in \Nz} [ J^n(A),  \cM_\infty B^{\dS^n}]. \]
Here, the connecting maps in the colimit are given by the assignment described above.
Note that $\cM_\infty B^{\dS^n}$ is an $\ind-\hal$-algebra indexed over $\Nz^{n+1}$. A priori it is possible to write
the definition of $kk$ without any mentioning of $\ind$-algebras, but it would involve too many colimits to be understandable at
an intuitive level.
\mn
Note that the sets appearing in the colimit carry abelian group structures and that the maps in the diagram of the colimit
are homomorphisms of abelian groups. It is also easy to see that the semigroup product coming from direct sum
of matrix blocks agrees with the group structure described above. This can be proved by giving an explicit rotational
homotopy. However, we omit the proof, since it also follows immediately from Proposition \ref{basic} once the
split-exactness (see Corollary \ref{split-exact}) is proved. Thus, $kk(A,B)$ is an abelian group.

\begin{remark}
The definition of $kk$-theory we have just given is relative to the ground ring $\hal$, and to our fixed choices of $\aha$ and of the underlying
category $\cU$. However it turns out that in the cases when we are able to compute them, the resulting $kk$-groups are independent of these choices. For example,
we show in \ref{proof_main} that for commutative $\hal$ and $\aha=\caha$, the groups $kk(\hal,A)$ do not depend on either $\hal$
or $\cU$. For another example, assume $\hal$ is commutative, fix $\cU=\hal\bimod$, and write $kk^b$ and $kk^{bc}$ for the a priori different
$kk$-groups which result from the choices $\aha=\gaha$ and $\aha=\caha$. If $A\in \caha$, then the value of $T(A)$ is independent of whether we choose $\aha=\caha$ or
$\aha=\gaha$, and, moreover, $T(A)\in \caha$. It follows that for all $n$, the algebra $J^n(A)$ is central, and does not depend on the choice of $\aha$.
On the other hand, the inclusion $\caha\subset\gaha$ has a right adjoint given by
\[
Z_\hal(B)=\{b\in B:[\hal,b]=0\}
\]
Moreover, one checks that $Z_\hal(B^{\dS^n})=Z_\hal(B)^{\dS^n}$. It follows that
\begin{equation}\label{noncentral}
kk^b(A,B)=kk^{bc}(A,Z_\hal B)\qquad (A\in\caha,\ \ B\in\gaha)
\end{equation}
Note that this can be further generalized to noncommutative $\hal$, by considering the commutativization $\hal_{ab}=\hal/\langle [\hal,\hal]\rangle$; details
are left to the reader.
\end{remark}

\subsection{Composition product}\label{subsec:compopro}

We write $\gamma_A: J(A^{\dS^1}) \to J(A)^{\dS^1}$ for the map discussed in Corollary \ref{natural_map}.

\begin{lemma}\label{left-right}
The map $\rho$ of \ref{subsec:loop} induces a natural transformation $J(?)\to(?)^{\dS^1}$ of endofunctors of $\ac$. That is, if
$f:A\to B\in \ac$, then
\[
f^{\dS^1} \rho_A=\rho_B J(f)\in\{J(A),B^{\dS^1}\}.
\]
\end{lemma}

\begin{lemma}\label{natural}
Let $A$ be an $\ind-\hal$-algebra. Then
\[ \gamma_A  J(\rho_{A})=\rho_{JA}: J^2(A) \to J(A)^{\dS^1}. \]
\end{lemma}

\begin{theorem}\label{kkcat}
Let $A,B$ and $C$ be $\ind-\hal$-algebras. There is an associative composition product
\[kk(B,C) \times kk(A,B) \to kk(A,C),\] which extends the composition of algebra homomorphisms.
\end{theorem}

\begin{proof}
Let $[\alpha] \in kk(B,C)$ be represented by $\alpha: J^n(B) \to C^{\dS^n}$ in $\acm$. Similarly,
let $[\beta] \in kk(A,B)$ be represented by $\beta: J^m(A) \to B^{\dS^m} $. We define the product
to be represented by the composition
\[ \xymatrix{J^{n+m}(A) \ar[r]^{J^n(\beta)} & J^n( B^{\dS^m}) \ar[r] & J^n(B)^{\dS^m} \ar[r]^{ \alpha^{\dS^m}}
& C^{\dS^{n+m}}. }\]
The fact that this assignment is well-defined in the colimit and also the
associativity of the product operation
is a formal consequence of Lemma \ref{left-right}, Lemma \ref{natural} and the naturality of the map
$\gamma_A: J(A^{\dS^1}) \to J(A)^{\dS^1}$ discussed in Corollary \ref{natural_map}.
\end{proof}

\begin{definition}
We write $kk$ for the category whose objects are those of $\aha$, and where the homomorphisms from $A\to B$ are given
by
\[
\hom_{kk}(A,B)=kk(A,B).
\]
\end{definition}

Let us note that there is a sequence of maps
\[ \xymatrix{ \hom_{\aha}(A,B) \ar[r]& [A,B] \ar[r] & \{A,B\} \ar[r]& kk(A,B) }\] which extend to functors between the
different categories. In particular the composite defines a functor $j:\aha\to kk$. To alleviate notation, and since
$j$ is the identity on objects, we write $A$ for $j(A)$. Moreover, when no confusion can arise, we also write $f$
for the image $j(f)\in kk(A,B)$ of a homomorphism $f\in\hom_{\aha}(A,B)$.
\mn
A priori, $kk$-theory is only defined for \hal-algebras. However, if $(A,J)$ is an $\ind-\hal$-algebra for which all
structure maps are $kk$-equivalences, we can equally well speak about its $kk$-groups. In some cases we will write $kk$-groups
with an $\ind-\hal$-algebra in the right argument. These are defined via the colimit of the induced diagram of $kk$-groups, i.e.
\[ kk(A,(B,J)) \deq \colim_{j \in J} kk(A,B_j).\]

The statement that two $\ind-\hal$-algebras $B$ and $C$ are $kk$-equivalent, has the rather strict meaning, that
all structure maps of $B$ and $C$ are $kk$-equivalences and all morphisms that constitute the morphism of
$\ind-\hal$-algebras are $kk$-equivalences as well.

\subsection{Excision}\label{subsec:exci}

\begin{lemma} \label{exact}
Let $f: A \to B$ be a morphism in $\aha$ and let $C$ be an \hal-algebra.
The sequence \[ \xymatrix{ kk(C,\cP_f) \ar[r]^{j(\ev_1)\cdot } & kk(C,A) \ar[r]^{j(f)\cdot} & kk(C,B) } \] is exact in the middle.
\end{lemma}
\begin{proof} The statement is a rewriting of the definition of being null-homotopic. \end{proof}

Let $f: A \to B$ be an $F$-split surjection of \hal-algebras. The $\ind-\hal$-algebra $\cP_f$  is indexed over $\Nz$ and
$\cP_{f,n} = B^{\sd^n \Delta^1} \oplus_B A$. The natural inclusions $\iota_n: \ker(f) \to P_{f,n}$ assemble to
a morphism of \hal-algebras $\iota: \ker(f) \to \cP_f$.

\begin{lemma}\label{ident}
Let $f: A \to B$ be an $F$-split surjection. Then the inclusion $\iota_n: \ker(f) \to \cP_{f,n}$ is invertible in $kk$
($n \geq 1$).
\end{lemma}
\begin{proof} Let $n \in \Nz$ be fixed.
There is a commutative diagram of extensions as follows.
\[ \xymatrix{ {\ker}(f)^{\sd^n S^1} \ar[r] \ar[d] & {\cP}_n \ker(f) \ar[r] \ar[d] & {\ker}(f) \ar[d]^{\iota} \\
{\ker}(f)^{\sd^n S^1} \ar[r] \ar[d]^{\iota^{\sd^n S^1}} & {\cP}_n A \ar[r] \ar[d] & {\cP}_{f,n} \ar[d]^{\id_{\cP_{f,n}}} \\
{\cP_{f,n}}^{\sd^n S^1}    \ar[r] & {\cP}_n \cP_{f,n} \ar[r] & {\cP}_{f,n}} \]

Here, an explicit homomorphism $\cP A  \to \cP \cP_f$ is needed. It is induced from a natural transformation $(?)^{\sd^\bullet \Delta^1} \to
((?)^{\sd^\bullet \Delta^1})^{\sd^\bullet \Delta^1}$. The required natural transformation has a level presentation
$(?)^{\sd^n \Delta^1} \to ((?)^{\sd^n \Delta^1})^{\sd^n \Delta^1}$ which we will describe. By Proposition \ref{power=tensor},
to give a natural map $A^{\sd^n \Delta^1} \to (A^{\sd^n \Delta^1})^{\sd^n \Delta^1}$ it is sufficient to
give $n^2$ maps $\eta_{k,l}: \Zz^{\sd^n \Delta^1} \subset \oplus_{i=1}^n \Zz[s] \to \Zz[t_1,t_2]$ for $1 \leq k,l \leq n$
which paste together in the correct way. We map the $i$-th summand by
\[ \eta_{k,l}(s) = \left\{ \begin{array}{cl}  t_1 & l<i=k \\
                            1- t_1 -t_2 +t_1t_2 & i= l = k \\
                            t_2 & k<i=l. \end{array} \right\}  \]
One checks that these indeed paste together correctly.
\mn
By Proposition \ref{classmap}, the following diagrams of classifying maps commute up to homotopy.

\[ \xymatrix{ J(\ker(f)) \ar[d]^{J(\iota)} \ar[r]^{\rho_{\ker(f)}} &{\ker (f)}^{\sd^n S^1} \ar[d]^{\id_{\ker(f)^{\sd^n S^1}}}
&& J(\cP_{f,n}) \ar[r]^\alpha \ar[d]^{\id_{J(\cP_{f,n})}} & {\ker}(f)^{\sd^n S^1} \ar[d]^{\iota^{\sd^n S^1}} \\
J(\cP_{f,n}) \ar[r]^\alpha & {\ker}(f)^{\sd^n S^1} && J(\cP_{f,n}) \ar[r]^{\rho_{\cP_{f,n}}} & {\cP_{f,n}}^{\sd^n S^1}}\]

This together implies that the class $[\alpha] \in kk(\cP_{f,n},\ker(f))$ is an inverse of $j(\iota_n)$.
\end{proof}

It follows from the lemma above that all maps in the inductive system defining $\cP_f$ are $kk$-equivalences.
Note also, that as far as $kk$-theory is concerned, there is no need
to distinguish between $\cP_f$ and $P_f$, as long as $f$ is surjective. We will later see that the inclusion $P_f \ton{\iota} \cP_f$ is
a $kk$-equivalence in general.
\mn
Applying the Lemma to the $F$-split surjection $\cP_n A \ton{\ev^n_1} A$, we get that
the inclusion $A^{\sd^n S^1} \to \cP_{\ev^n_1}$ is a $kk$-equivalence. Note that $\cP_{\ev^n_1}=\cP_{\ev_1,n}$, so that, again by the Lemma,
we get that the inductive system
\[ \sd^{\bullet}A: \Omega A \to A^{\sd^1 S^1}\to A^{\sd^2 S^1}\to A^{\sd^3 S^1} \to \cdots\]
consists of $kk$-equivalences. In particular $\iota: \Omega A \to A^{\dS^1}$ is a kk-equivalence.

\begin{corollary}\label{path-algebras} Let $f: A \to B$ be a morphism of \hal-algebras, $\pi_f:\cP_f\to A$ the canonical map.
The inclusion of $\ind-\hal$-algebras $B^{\dS^1} \to \cP_{\pi_f}$ is invertible in $kk$.
\end{corollary}
\begin{proof} The proof is an application of Lemma \ref{ident} in the case of the surjection $B \to 0$ and
$\ev_1: \cP_{f,n} \to A$.
\end{proof}

\begin{corollary} \label{split-exact}
Let $D$ be an \hal-algebra. The functor $kk(D,?)$ is split-exact. \end{corollary}
\begin{proof}
Let \[
\xymatrix{ A \ar[r] & B \ar[r]^g & C}
\]
be a split extension. By Lemma \ref{exact},
we know that
\[
 \xymatrix{0 \ar[r] & kk(D,P_g) \ar[r] & kk(D,B) \ar[r] & kk(D,C) \ar[r] & 0}
\]
is exact in the middle.
Since the surjection $g: B \to C$ is split, it is also exact at $kk(D,C)$. By Lemma \ref{ident},
the natural inclusion $\iota: A \to P_g$ is an equivalence in $kk$, so that $kk(D,P_g)$ identifies with
$kk(D,A)$. It remains to show that
\[
\xymatrix{kk(D,P_g) \ar[r]^{j(\pi_g)\cdot}& kk(D,B)}
\]
is injective.
By Lemma \ref{exact}
\[
 \xymatrix{ kk(D,P_{\pi_g}) \ar[r] & kk(D,P_g) \ar[r] & kk(D,B) }
\]
is exact in the middle. Furthermore,
by Corollary \ref{path-algebras}, the inclusion $\Omega C \to P_{\pi_g}$ is a $kk$-equivalence. Now, the composition
$\Omega C \to P_g$ factors through $PC$, since $g$ is split. This shows that
\[
 \xymatrix{ kk(D,P_{\pi_g}) \ar[r] & kk(D,P_g)}
\]
is zero and thus, by exactness, $kk(D,P_g) \to kk(D,B)$ has to be injective. This finishes the proof.
\end{proof}

\begin{corollary} \label{iterated_cone}
Let $f: A \to B$ be an $F$-split surjection. Consider the natural inclusion $\iota_f: \Omega B \to P_f$.
The natural map $\pi: P_{\iota_f} \to \Omega A$ is a $kk$-equivalence.
\end{corollary}
\begin{proof}
One shows that $\pi$ is a split-surjection and that the kernel is contractible. Then, the result follows from Corollary \ref{split-exact}.
\end{proof}

\begin{theorem}\label{theo:exci_1variable}
Let
\[
\xymatrix{ A \ar[r]^f & B \ar[r]^g & C}
\]
be an $F$-split extension. Then the following sequence is exact.
\[
\xymatrix{ kk(D, \Omega B) \ar[r]^{j(\Omega g)\cdot } & kk (D,\Omega C) \ar[r]^\partial & kk(D,A) \ar[r]^{j(f)\cdot}
& kk(D,B) \ar[r]^{j(g)\cdot} & kk(D,C) }
\]
Here, the map $\partial: kk(D,\Omega C) \to kk(D,A)$ is given by the composition of the first map in the diagram below followed
by the inverse of the second
\[
 \xymatrix{ kk(D,\Omega C) \ar[r] & kk(D,P_g) & kk(D,A) \ar[l]_{\sim}}.
\]
\end{theorem}
\begin{proof} This follows immediately from the iteration of Lemma \ref{exact} and the identifications which are possible due
to Corollaries \ref{path-algebras} and \ref{iterated_cone}. \end{proof}

\begin{theorem} \label{sequence}
Let
\[
\xymatrix{ A \ar[r]^f & B \ar[r]^g & C}
\]
by an $F$-split extension. The following sequence is natural and exact.
\[
\xymatrix{ kk(C,D) \ar[r]^{\cdot j(g) } & kk (B,D) \ar[r]^{\cdot j(f)} & kk(A,D) \ar[r]^{\partial} & kk(\Omega C,D)
\ar[r]^{\cdot j(\Omega g)} & kk(\Omega B,D) }
\]
Here, the map $\partial: kk(A,D) \to kk(\Omega C,D)$ is given by the composition of the inverse of the first map below followed
by the second map
\[
\xymatrix{ kk(A,D) & kk(P_g,D) \ar[l]_{\sim} \ar[r]& kk(\Omega C,A)}.
\]
\end{theorem}
\begin{proof}
We start by showing that for any $F$-split surjection $f: A \to B$, the sequence
\[\xymatrix{kk(B,D) \ar[r]^{\cdot j(f)} & kk(A,D) \ar[r]^{\cdot j(\pi_f)}& kk(P_f,D)}\]
is exact in the middle.
Having proved this, the result follows from Corollaries \ref{path-algebras} and \ref{iterated_cone}.
\mn
Let $[\alpha]$ be an element in $kk(A,D)$ which vanishes after precomposing with $j(\pi_f)$. Let $[\alpha]$
be represented by $\alpha: J^n(A) \to \cM_\infty D^{\dS^n}$. We can choose $n$ so large, that $\alpha  J^n(\pi_f)$
is null-homotopic. This leads to the existence of a commutative diagram of extensions as follows. (Note here, that the functor
$J$ preserves $F$-split extensions.)

\[ \xymatrix{{\ker}(J^n(\pi_f)) \ar[r] \ar[d] & J^n(P_f) \ar[r]^{J^n(\pi_f)} \ar[d]& J^n(A)\ar[d]^{\alpha} \\
{\cM}_\infty D^{\dS^{n+1}} \ar[r] & {\cM}_\infty \cP D^{\dS^n} \ar[r] &{\cM}_\infty D^{\dS^n}  }\]

Consider the composite
\[
\beta:J^{n+1}(B) \to J^n(\Omega B) \to \ker(J^n(\pi_f)) \to \cM_\infty D^{\dS^{n+1}}.
\]
Write $[\beta]$ for its class in $kk(B,D)$. We claim that $j(f)  [\beta]=[\alpha]$. This follows
by diagram chasing, using the uniqueness of the classifying map of an extension up to homotopy.

 \end{proof}

\begin{lemma} \label{isom-stab}
Let $A$ and $B$ be \hal-algebras. There is a natural isomorphism
\[ \xymatrix{\Omega: kk(A,B) \ar[r]^{\sim} & kk(\Omega A, \Omega B).}\]
\end{lemma}
\begin{proof} Let $[\alpha]$ be a class in $kk(A,B)$ which is represented by $\alpha: J^n(A) \to B^{\dS^n}$. To it we associate
the class of the composition $J^n(A^{\dS^1}) \to J^n(A)^{\dS^1} \ton{\alpha^{\dS^1}} B^{\dS^{n+1}}$ in
$kk(A^{\dS^1},B^{\dS^1}) \cong kk(\Omega A, \Omega B)$. We call this assignment \textit{looping}.
Because the following diagram commutes in $\ac$, looping is well-defined at the level of $kk$-theory.
\[
\xymatrix{J^{n+1}(A^{\dS^1}) \ar[d]^1 \ar[r] & J(J^n(A)^{\dS^1}) \ar[d]^{\gamma_{J^n(A)}} \ar[r]^{J(f^{\dS^1})}
& J(B^{\dS^{n+1}}) \ar[d]^{\gamma_{B^{\dS^{n}}}} \ar[rr]^{\rho_{B^{\dS^{n+1}}}} && B^{\dS^{n+2}}  \ar[d]^1 \\
J^{n+1}(A^{\dS^1}) \ar[r] & J^{n+1}(A)^{\dS^1} \ar[r]^{J(f)^{\dS^1}}  & J(B^{\dS^{n}})^{\dS^1}
\ar[rr]^{\rho_{B^{\dS^{n}}}^{\dS^1}} && B^{\dS^{n+2}}}.
\]
Indeed, the lower row represents the effect of looping after having stabilized. The upper row represents the
stabilization in $kk$ after looping.

It remains to show that looping induces an isomorphism. Given a class in $kk(A^{\dS^1},B^{\dS^1})$ which is represented
by $\beta: J^n(A^{\dS^1}) \to B^{\dS^{1+n}}$, we assign to it the class in $kk(A,B)$ of the composition
$$J^{n+1}(A) \ton{J^n(\rho_A)} J^n(A^{\dS^1}) \ton{\beta} B^{\dS^{1+n}}.$$ Again, this
is well defined at the level of $kk$. We call this procedure \textit{delooping}. The following diagram shows that delooping
after having looped gives stabilization, and hence the identity in $kk$.

\[\xymatrix{J^n(A^{\dS^1})  \ar[r]^{\gamma} & J^n(A)^{\dS^1} \ar[r]^{f^{\dS^1}} & B^{\dS^{n+1}} \\
J^{n+1}(A) \ar[u]^{J^n(\rho_A)} \ar[r]^1 & J^{n+1}(A) \ar[u]^{\rho_{J^n(A)}} \ar[r]^{J(f)} &
J(B^{\dS^n}) \ar[u]^{\rho_{B^{\dS^n}}} }\]

Indeed, starting from the lower left corner and going up and right gives the delooping of the loops, going right
and up gives stabilization in $kk$. The squares commute up to homotopy by Lemmas \ref{natural} and \ref{left-right}.
It is similarly shown that delooping first and then looping  also gives stabilization. This finishes the proof.
\end{proof}

\begin{lemma}\label{omegaisfunctor} Looping is compatible with composition. Thus it defines an endofunctor of the category $kk$, which sends
an algebra $A$ to $\Omega A$. We use the notation $\Omega: kk \to kk$.\end{lemma}

Lemma \ref{isom-stab} shows that $\Omega$ is fully faithful. We will see later that $\Omega$ is also essentially surjective, i.e. an equivalence of
categories.

\begin{lemma}\label{lemma_about_rho}
Let $A$ be an \hal-algebra. The natural map $\rho_A: J(A) \to \Omega A$ induces a $kk$-equivalence.
\end{lemma}
\begin{proof}
Note that $\rho_A$ is part of a natural map between extensions, the universal extension and the path extension. In both extensions,
the middle term is contractible and hence isomorphic to zero in $kk$. The desired result follows from the naturality of
the exact sequences associated to an $F$-split extension, which were established in Theorem \ref{sequence}.
\end{proof}

The following lemma provides a useful description of an element in $kk(A,B)$ which we will use it in the subsequent sections.

\begin{lemma} \label{comp_loops}
Let $[f]$ be an element in $kk(A,B)$, which is represented by $f: J^n(A) \to B^{\dS^n}$. The composite of the first two maps in the
following diagram, followed by the inverse of the third, induces $\Omega^n [f]$ in $kk$-theory.
\[
 \xymatrix{\Omega^n A \ar[rr]^{j(\rho^n_A)^{-1}} && J^n(A) \ar[r]^f  & B^{\dS^n}& \Omega^n B \ar[l]_{\sim} }.
\]

\end{lemma}

\mn

\subsection{The pseudo-inverse to $\Omega$}\label{subsec:deloop}

We claimed that the functor $\Omega$ is essentially surjective. Using the following theorem, we will provide a proof of this.

\begin{theorem}\label{theo:infsum_zero} Let $(A,\phi^{\infty})$ be a infinite sum $\hal$-algebra, and $B\vartriangleleft A$ an ideal such that $\phi^{\infty}(B)\subset B$.
Then $B$ is $kk$-equivalent to zero. \end{theorem}

\begin{proof} From Lemma \ref{wagoner_lemma} and Proposition \ref{prop:multiplier}, we get that $j(\id_B \oplus \phi^{\infty}_{|B}) = j(\id_B) +
j(\phi^{\infty}_{|B}) \in kk(B,B).$
On the other hand, by definition of infinite sum \hal-algebra, $j(\id_B \oplus \phi^{\infty}_{|B}) = j(\phi^{\infty}_{|B})$, so that
$j(\id_B)$ has to vanish. This finishes the proof.
\end{proof}

\begin{corollary} \label{inverse}
Let $A$ be an \hal-algebra. Then $\Sigma A $ is a delooping of $A$, i.e. there is a natural $kk$-equivalence
\[ \xymatrix{{\Omega} \Sigma(A) \ar[r]^{\sim} &  A.} \]
\end{corollary}
\begin{proof} Immediate from \ref{infinitesum}, its proof, and the theorem above.
\end{proof}

\mn

We use the stabilization isomorphism from Lemma \ref{isom-stab} and Corollary \ref{inverse},
to define $\Zz$-graded $kk$-groups as follows:

\[ kk_n(A,B) \deq \left\{  \begin{array}{ll} kk(A,\Omega^n B) & \mbox{for} \, n\geq 0 \\
kk(\Omega^{-n} A,B) \cong kk(A,\Sigma^{-n} B) & \,\mbox{for}\, n<0.  \end{array} \right\}
\]
\sn
In order to organize the properties of $kk$, we state that $kk_*$ defines a bivariant homology theory in the
sense of \cite{thom-thesis}, i.e. it produces in each variable long exact sequences out of $F$-split extensions.
In Subsection \ref{subsec:kktria} below, we will show that $kk$ is indeed a triangulated category with a suitable class of distinguished triangles. This puts all
the homological structure in an appropriate framework.

\mn

\subsection{$kk$ as a triangulated category}\label{subsec:kktria}
For the definition of triangulated category we refer to \cite{Keller}. The original definition goes back to D. Puppe (\cite{puppe})
and J.-L. Verdier (\cite{verdier}),
who studied the main motivating examples of triangulated categories, i.e. the stable homotopy category and derived categories.
Over the last few years, several papers have appeared, in which triangulated categories of rings or algebras (in particular operator algebras) are
considered (e.g. \cite{puschnigg, meyer, thom-thesis}). Here, we will show that $kk$ is a triangulated category.
\mn
The definition of triangulated category we shall be using is that of \cite{Keller}; axioms shall be recalled during the course of the proof of Theorem \ref{triangcat}.
This definition is a slight variant of that given for example in A. Neeman's book \cite{neeman}; instead of requiring that the endofunctor $\Sigma$ be invertible, we just require that it be an equivalence of categories. It has
been shown (see for example in the Appendix of \cite{margolis}) that this yields no greater generality, i.e. each triangulated category in this weak sense is triangulated equivalent
to an ordinary triangulated category.
\mn
The definition of triangulated category is self-dual, i.e. the opposite category
of a triangulated category inherits a natural triangulated structure. It is therefore enough ot check that $kk^{op}$ is triangulated, which is
much more natural in our setting. In turn, checking this amounts to showing that $kk$ satisfies the axioms opposite to the usual axioms.
These are most naturally stated in terms of a pseudo-inverse $\Omega$ of $\Sigma$.
\mn
\begin{definition} \label{triang}
A diagram \[\xymatrix{{\Omega} C \ar[r]& A \ar[r] & B \ar[r] & C }\] of morphisms in $kk$ is called a distinguished triangle,
if it isomorphic in $kk$ to the path sequence
\[\xymatrix{ {\Omega} B' \ar[r]^{j(\iota)}& P_f \ar[r]^{j(\pi_f)} & A' \ar[r]^{j(f)} & B' }\]
associated with a homomorphism $f: A' \to B'$ of \hal-algebras.
\end{definition}

\begin{theorem} \label{triangcat}The category $kk$ is triangulated with respect to the endofunctor $\Omega: kk \to kk$ and the class
of distinguished triangles specified in Definition \ref{triang}.
\end{theorem}
\begin{proof} First of all, we need to show that $\Omega$ is an equivalence. This follows from Lemma \ref{isom-stab} and Corollary \ref{inverse}.
Next we verify that all the requisite axioms for a triangulated category are satisfied.

\begin{axiom}[\bf TR0]
Any distinguished triangle which is isomorphic to a distinguished triangle is a distinguished triangle. The
triangle
\[\xymatrix{{\Omega} A \ar[r]& 0 \ar[r] & A \ar[r]^{1} & A }\]
is distinguished, for every object $A$ in $kk$. \end{axiom}
The first assertion is clear from Definition \ref{triang}. To prove the second, consider the path sequence of
$1_X$ and note that $P_{1_X}$ is isomorphic to $0$ in $kk$.
\begin{axiom}[{\bf TR1}]
For any morphism $\alpha: A \to B$ in $kk$, there exists a distinguished triangle of the form
\[\xymatrix{{\Omega} B \ar[r]& C \ar[r] & A \ar[r]^\alpha & B. }\]
\end{axiom}
Let $\alpha$ be represented by the $\ind-\hal$-algebra homomorphism $f: J^n(A) \to B^{\dS^n}$. By Lemma \ref{comp_loops}
the following diagram commutes
\[\xymatrix{   {\Omega}^n A \ar[r]^{\Omega^n \alpha}  & {\Omega}^n B \ar[d] \\
J^n(A) \ar[r]^f \ar[u]^{\rho^n_A} & B^{\dS^n}}\]
in $kk$-theory. This implies that $j(\Sigma^n(f))$ is isomorphic to $\alpha$ in the arrow category of $kk$. The path sequence of
$\Sigma^{n}(f)$ provides a sequence which contains $\alpha$ up to isomorphism as desired.
\begin{axiom}[{\bf TR2}]
Consider the two triangles
\[\xymatrix{{\Omega} C \ar[r]^f& A \ar[r]^g& B \ar[r]^h & C \\
{\Omega} B \ar[r]^{-\Omega(h)} & {\Omega} C \ar[r]^{-f} & A \ar[r]^{-g} & B. }\]
If one is distinguished, then so is the other.
\end{axiom}
We call the lower triangle the rotation of the upper triangle. It is obvious that the threefold rotation of the path-sequence
of $f$ is isomorphic to the path-extension of $\Omega f$. Also, the path-sequence of $\Omega \Sigma f$ is isomorphic
to the path-extension of $f$. This implies that the threefold rotation of a triangle is distinguished if and only if the triangle is
distinguished.
Thus, it suffices to show that the rotation of a distinguished triangle is distinguished. \mn
By Definition \ref{triang}, we may assume that the first triangle is equal to
\[\xymatrix{ {\Omega} B' \ar[r]^{j(\iota)}& P_f \ar[r]^{j(\pi_f)} & A' \ar[r]^{j(f)} & B'. }\]
The second triangle is then isomorphic to
\[\xymatrix{{\Omega} A' \ar[r]^{-j(\Omega(f))} & {\Omega} B' \ar[r]^{j(\iota)}& P_f \ar[r]^{j(\pi_f)} & A'. }\]
By Lemma \ref{path-algebras} the natural map of \hal-algebras $\Omega B'  \to P_{\pi_f}$
is a $kk$-equivalence and makes the diagram
\[\xymatrix{{\Omega}{A'} \ar[d]^1 \ar[r]^{-j(\Omega(f))} & {\Omega}{B'} \ar[d] \ar[r]^{j(\iota)}& P_f \ar[d]^1 \ar[r]^{j(\pi_f)} & A' \ar[d]^1\\
{\Omega}{A'} \ar[r]^{\iota} & P_{\pi_f}\ar[r] & P_f \ar[r]^{j(\pi_f)} & A'}\]
commute. This finishes the proof of axiom TR2.

\begin{axiom}[{\bf TR3}]
For any commutative diagram
\[\xymatrix{{\Omega} C \ar[d]^{\Omega k}\ar[r]^f& A \ar[r]^g& B \ar[r]^h \ar[d]^l & C \ar[d]^{k} \\
{\Omega} C' \ar[r]^{f'}& A' \ar[r]^{g'}& B' \ar[r]^{h'} & C', } \] there exists a
filler $j: A \to A'$ which makes the whole diagram commutative.
\end{axiom}
We may assume that the horizontal sequences are path sequences. Furthermore, choosing $n$ sufficiently large, we may pick representatives
of $l$ and $k$ of the form $a: J^n(B) \to \Omega^n B'$ and $b: J^n(C) \to \Omega^n C'$. Increasing $n$ even more if necessary,
we may assume that the diagram of $\hal$-algebras
\[\xymatrix{ J^n(B) \ar[r]^{J^n(h)} \ar[d]^a & J^n(C) \ar[d]^b \\
{\Omega}^n B' \ar[r]^{\Omega^n h'} & {\Omega}^n C' }\] commutes up to homotopy.
It follows from the properties of $\Sigma$ and Lemma \ref{comp_loops}, that
\[\xymatrix{ \Sigma^n J^n(B) \ar[rr]^{\Sigma^nJ^n(h)} \ar[d]^{\Sigma^n a} && {\Sigma^n}J^n(C) \ar[d]^{\Sigma^nb}\\
{\Sigma^n}\Omega^n B' \ar[rr]^{\Sigma^n\Omega^n h'} && {\Sigma^n}\Omega^n C' }\] also
commutes up to homotopy and is isomorphic in $kk$ to the right-hand square in the diagram of triangles. It follows that the path sequences
of $\Sigma^nJ^n(h)$ and $\Sigma^n\Omega^n h'$ are isomorphic to those of $h$ and $h'$. To finish the proof, choose a null-homotopy for
$P_{\Sigma^nJ^n(h)} \to \Sigma^n\Omega^n C'$ and use it to construct a map
\[
P_{\Sigma^nJ^n(h)} \to P_{\Sigma^n\Omega^n h'}= P(\Sigma^n\Omega^n B') \oplus_{\Sigma^n\Omega^n C'}
\Sigma^n\Omega^n B'.
\]

\begin{axiom}[{\bf TR4}]
Let $\alpha: A \to B$ and $\beta: B \to C$ be morphisms in $kk$. We set $\gamma= \beta  \alpha$.
There exists a commuting diagram
\[\xymatrix{
{\Omega}^2 C \ar[r]\ar[d]&{\Omega} D'\ar[r]\ar[d]& {\Omega} B\ar[r]^{\Omega \beta}\ar[d]& {\Omega} C \ar[d]\\
0\ar[r] \ar[d]&D'''\ar[r]^1 \ar[d]^h &D'''\ar[r] \ar[d]&0\ar[d] \\
{\Omega} C\ar[r]^j \ar[d]^1& D''\ar[r] \ar[d]&A \ar[d]^{\alpha} \ar[r]^{\gamma}& C \ar[d]^1\\
{\Omega} C\ar[r]& D'\ar[r]&B \ar[r]^{\beta }& C\\       }\]
in which each row and column is an exact triangle. Furthermore, the square
\[\xymatrix{
{\Omega} B \ar[r]^{\Omega \beta} \ar[d]& {\Omega} C \ar[d]^j \\
D''' \ar[r]^h & D''  }\] commutes in $kk$.
\end{axiom}
Like in the proof of TR3, we replace the commuting square in the lower right corner by an isomorphic square which has a lift to a
square in the category of \hal-algebras which commutes up to homotopy. Using a path-algebra, we may require that
it strictly commutes.
\mn
The lower horizontal sequences are defined to be the path sequences of $\beta$ and $\gamma$, i.e. $D'=P_{\beta}$ and $D''=P_{\gamma}$.
Let $\eta: P_{\gamma} \to P_{\beta}$ be the induced map. The square
\[\xymatrix{ P_{\gamma} \ar[r] \ar[d]^{\eta}& A \ar[d]^{\alpha} \\
P_{\beta} \ar[r] & B  }\] commutes, so that there is a natural map between the path sequences of $\eta$ and $\alpha$. We put these
into the middle vertical sequences. Note that the induced map $P_{\eta} \to P_{\alpha}$ is a split surjection with contractible kernel, hence
a $kk$-equivalence, so that we can put $D'''=P_{\eta}$. It is now clear
that all rows and columns are exact triangles and all squares in the big diagram commute.
\mn
The commutativity of the small diagram follows from the commutativity of
\[\xymatrix{ {\Omega} B\ar[r]^{\Omega \beta} \ar[d] & {\Omega} C \ar[d] \\
P_{\alpha} \ar[r] & P_{\gamma}  }\] and the identification $P_{\eta} \to P_{\alpha}$. This finishes the proof of Axiom TR4.
\end{proof}

\subsection{Universal properties}\label{subsec:uniprop}
Let $\cE$ be the class of all $F$-split extensions
\begin{equation}\label{extensionE}
(E):A\to B\to C
\end{equation}
in $\aha$, and $\cT=(\cT,\Omega)$ a triangulated category.
An {\it excisive homology theory}
for $\hal$-algebras with values in $\cT$ consists of a functor
$X:\aha\to \cT$, together with a collection
$\{\partial_E:E\in\cE\}$ of maps $\partial_E^X=\partial_E\in\hom_{\cT}(\Omega
X(C), X(A))$. The maps $\partial_E$ are to satisfy the following requirements.
\sn
\noindent{i)} For all $E\in \cE$ as above,
\[
\xymatrix{\Omega
X(C)\ar[r]^{\partial_E}&X(A)\ar[r]^{X(f)}&X(B)\ar[r]^{X(g)}& X(C)}
\]

is a distinguished triangle in $\cT$.
\sn \noindent{ii)} If
\[
\xymatrix{(E): &A\ar[r]^f\ar[d]_\alpha& B\ar[r]^g\ar[d]_\beta& C\ar[d]_\gamma\\
          (E'):&A'\ar[r]^{f'}& B'\ar[r]^{g'}& C'}
\]
is a map of extensions, then the following diagram commutes
\[
\xymatrix{{\Omega} X(C)\ar[d]_{{\Omega} X(\gamma)}\ar[r]^{\partial_E}& X(A)\ar[d]^{X(\alpha)}\\
 \Omega X(C')\ar[r]_{\partial_{E'}}& X(A)}
\]
We say that the functor $X:\aha\to\cT$ is {\it homotopy invariant}
if it maps homotopic homomomorphisms to equal maps, or
equivalently, if for every $A\in\aha$, $X$ maps the inclusion
$A\subset A[t]$ to an isomorphism. Call $X$ {\it $M_\infty$-stable} if for every $A\in\aha$,
it maps the inclusion $\iota_\infty:A\to M_\infty A$ to an isomorphism. Note that if $X$ is
$M_\infty$-stable, and $n\ge 1$, then $X$ maps the inclusion $\iota_n:A\to M_nA$ to an isomorphism.

\begin{example}
Let $E$ be the $F$-split extension \eqref{extensionE}, and $c_E\in kk(JC,A)$ the class of the classifying map.
Define
\begin{equation}\label{collect}
\partial_E\deq c_E\circ \rho_A^{-1}\in kk(\Omega C,A).
\end{equation}
Then the canonical functor $j:\aha\to kk$,
together with the collection $\{\partial_E\}_{E\in\cE}$ is a homotopy invariant, $M_\infty$-stable, excisive homology
theory in the sense defined above.
\end{example}
The homotopy invariant, $M_\infty$-stable, excisive homology theories form a
category, where a homomorphism between the theories $X:\aha\to
\cT$ and $Y:\aha\to \cU$ is a triangulated functor $G:\cT\to\cU$ such
that
\begin{equation}\label{fobjects}
\xymatrix{{\aha}\ar[dr]_Y\ar[r]^X&{\cT}\ar[d]^G\\
                               &{\cU}}
\end{equation}
commutes, and such that for every extension \eqref{extensionE}, the natural isomorphism
$\phi:G(\Omega X(C))\to \Omega Y(C)$ makes the following into a commutative diagram
\begin{equation}\label{phi-partial}
\xymatrix{G(\Omega X(C))\ar[r]^(.6){G(\partial^X_E)}\ar[d]_\phi&Y(A)\\
             \Omega Y(C).\ar[ur]_{\partial^Y_E}& }
\end{equation}

\begin{theorem}\label{univ:triang} The canonical functor $j:\aha\to kk$, together with
the maps \eqref{collect}, is an initial object in the category of all excisive,
homotopy invariant and $M_\infty$-stable homology theories.
\end{theorem}
\begin{proof}
Let $X:\aha\to \cT$ be an excisive, homotopy invariant and $M_\infty$-stable theory. We have to show that there is a unique homomorphism
of theories $\bar{X}:j\to X$. By \eqref{fobjects}
we must have $\bar{X}(A)=X(A)$. Because $X$ is homotopy invariant, the connecting map associated to any $F$-split extension \eqref{extensionE}
with $B$ contractible must be an isomorphism. It follows from this that $X$ sends each subdivision map $\Omega^nB\to B^{\sd^pS^n}$ to an isomorphism.
Further,
by \eqref{phi-partial}, $\phi=\partial_l^{-1}$, the inverse of the connecting map associated to the
loop extension.  Let $\alpha\in kk(A,B)$ be represented by a homomorphism $J^nA\to B^{\sd^pS^n}\in\{J^nA,B^{\sd^pS^n}\}$. Then
$\Omega^n\alpha\in kk(\Omega^nA,\Omega^nB)$ is the class of the composite of $\rho^{-n}\in kk(\Omega^nA,J^nA)$ followed by $f$,
followed by the inverse $\mu$ of $\Omega^nB\to B^{\sd^pS^n}$. Because $X$ is homotopy invariant and $M_\infty$-stable, it extends uniquely to a functor on the
full subcategory of $\acm$ whose objects are the constant $\ind$-algebras; we shall abuse notation and write $X$ for this
extension. In particular $X$ is defined on $\mu^{-1} f\rho^{-n}$, and we must have $\bar{X}(\Omega^n(\alpha))=X(\mu^{-1} f\rho^{-n})$.
Moreover, the condition that $\bar{X}$ be triangulated determines that for the inverse $\Omega^{-1}$ of the bijection
\[
\Omega :\hom_\cT(X(A),X(B))\to \hom_\cT(\Omega X(A),\Omega X(B))
\]
we must have
\begin{equation}
\begin{split}\label{forcedef}
\bar{X}(\alpha)=&\Omega^{-n}(\phi^n \bar{X}(\Omega^n\alpha)\phi^{-n})\\
=&\Omega^{-n}(\partial_l^{-n} X(\mu)^{-1}X(f)X(\rho)^{-n}\partial_l^{n})
\end{split}
\end{equation}
Essentially the same argument as that of the proof of \ref{kkcat} shows that \eqref{forcedef} actually defines a functor
$\bar{X}:kk\to \cT$. It remains to show that the functor $\bar{X}$ is triangulated. In view of \ref{triang} and of the rotation axiom,
this boils down to proving that if $f\in\hom_\aha(A,B)$, then $X$ maps
\begin{equation}\label{path:seq}
\xymatrix{\Omega A\ar[r]^{\Omega f}&\Omega B\ar[r]^\iota& P_f\ar[r]^\pi &A}
\end{equation}
to a distinguished triangle in $\cT$. Consider the extension $E$ formed by $\iota$ and $\pi$; we have a
diagram
\[
\xymatrix{& \Omega X(B)\ar[d]^{\partial_l}& &\\
\Omega X(A)\ar[ur]^{\Omega X(f)}\ar[r]_{\partial_E}\ar[d]^{\partial_l}& X(\Omega B)\ar[r]_{X(\iota)}&X(P_f)\ar[r]^{X(\pi)}&X(A)\\
 X(\Omega A)\ar[ur]_{X(\Omega f)}&&&}
\]
In the diagram above the row is a distinguished triangle and the square on the left commutes, as do its upper and lower halves.
It follows from this and the axioms of a triangulated category, that $X$ applied to \eqref{path:seq} is a distinguished triangle.
\end{proof}

Next we give three examples which illustrate the applications of the theorem above.
\begin{example}
Algebraic $K$-theory is not excisive, nor is it $M_\infty$-stable or
homotopy invariant. However there is a variant defined by C. Weibel
\cite{weih}, called homotopy algebraic $K$-theory, and denoted
$KH$, which does have all these properties. The groups $KH_*(A)$
are defined as homotopy groups of a certain spectrum. There are
several homotopy equivalent choices for this spectrum; one of them
is recalled in \ref{subsec:khdef}. There is also a variant which yields a
functor $\kaha:\Ass\deq\Ass_\Zz\to \rm{Sp}^\Sigma$ with values in the
category of symmetric spectra (see \cite{thom-thesis}), and is
equipped with suitably compatible products
\[
\kaha(A)\land \kaha (B)\to \kaha(A\otimes B)
\]
In particular $\kaha$ maps $\hal$-algebras to
$\kaha(\hal)$-modules. Composing with the canonical localization
functor to the homotopy category $\cT$ of $\kaha(\hal)$-modules,
we get a functor $X:\aha\to\cT$ which satisfies the hypothesis of
Theorem \ref{univ:triang}. Thus we have a natural map
\begin{equation}\label{compa}
kk(A,B)\to KH(A,B)\deq [\kaha(A),\kaha(B)]_{\kaha(\hal)}\deq
\hom_\cT(X(A),X(B)).
\end{equation}
\end{example}

\begin{example}\label{chern}
Assume $\hal$ is a field of characteristic zero. Then bivariant
periodic cyclic cohomology (\cite{cq}), $(A,B)\mapsto HP^*(A,B)$,
is $2$-periodic, excisive, homotopy invariant and $M_\infty$-stable
in both variables. This bivariant theory is defined by means of a
functor $X^\infty$ from $\aha$ to pro-supercomplexes. The latter form
a closed model category (see \cite[4.3]{cv}) whose homotopy category $\cT$
is triangulated. One defines
\[
HP^*(A,B)\deq \hom_\cT(X^\infty(A),\Omega^*X^\infty(B)).
\]
The functor $X^\infty:\aha\to \cT$ satisfies the hypothesis of Theorem \ref{univ:triang},
(\cite{cv}) so that we have a product-preserving Chern character
\[
ch_*:kk_*(A,B)\to HP^*(A,B).
\]
\end{example}
\begin{example}
Assume $\hal$ is commutative, put $\aha=\caha$, and let $\cU$ be the category of $\hal$-modules (i.e. central bimodules). Then
$F$-split extensions remain $F$-split upon tensoring over $\hal$, so that for every algebra $B$ the functor
$\iota\circ(? \otimes_\hal B):\aha\to kk$ is excisive. On the other hand, it is clear that $\iota\circ(?\otimes_\hal B)$
is homotopy invariant and $M_\infty$-stable. Applying Theorem \ref{univ:triang} we obtain a functor $\overline{\otimes}_\hal B:kk\to kk$. It follows that
there is an associative product
\begin{equation}\label{kkext}
kk(A_1,B_1)\otimes kk(A_2,B_2)\to kk(A_1\otimes_\hal A_2,B_1\otimes_\hal B_2).
\end{equation}
\end{example}

\bn
After dealing with functors into triangulated categories, we now concentrate on homological (i.e. half-exact) functors with values in some abelian category.
Let $\Ac$ be an abelian category. A functor $G:\aha\to \Ac$ is called {\it half exact} if for every $F$-split extension \eqref{extension} the
sequence
\[
G(A)\to G(B)\to G(C)
\]
is exact.
\begin{theorem}\label{univ:hom} Let $\Ac$ be an abelian category, and $G:\aha\to\Ac$ a half exact, homotopy invariant, $M_\infty$-stable, additive functor.
Then there exists a unique homological functor $\bar{G}:kk\to \Ac$ such that the following diagram commutes.
\[
\xymatrix{{\aha}\ar[r]^j\ar[dr]_G&kk\ar@{.>}[d]^{\bar{G}}\\
                               &{\Ac}}
\]
\end{theorem}
\begin{proof}
Let $f:B\to C$ be the surjective map in \eqref{extensionE}, and $\iota:\Omega C\to P_f$ the inclusion. A standard argument shows (see for example \cite[\S21.4]{black})
that $G$ sends the natural map $l:A\to P_f$ to an isomorphism, and that for $\partial^G_E\deq G(l)^{-1}\circ G(\iota)$, the following
sequence is exact
\begin{equation}\label{H:path}
\xymatrix{G(\Omega B)\ar[r]^{G(\Omega f)}& G(\Omega C)\ar[r]^{\partial^G_E}& G(A)\ar[r]&G(B)\ar[r]^f&G(C)}
\end{equation}
Moreover one checks that the map $\partial^G_l$ corresponding to the loop extension is the identity map $G(\Omega C)\to G(\Omega C)$.
On the other hand, \eqref{H:path} implies that $\partial_E^G$ must be an isomorphism whenever $B$ is contractible. In particular this
applies to the connecting map $\partial^G_u$ associated to the universal extension. It also follows from \eqref{H:path} that
$G(\rho)$ is an isomorphism. Hence $\partial^G_E=G(c_E)G(\rho)^{-1}$ for every $F$-split extension \eqref{extensionE}.
By \eqref{collect}, we must define $\bar{G}(\partial_E)=\partial^G_E$. Now the argument of the proof of \ref{univ:triang} shows that
if $\alpha\in kk(A,B)$ then we have a unique way of defining $\bar{G}(\Omega^n\alpha)$.
Let $\tau:\Sigma\Omega\to \Omega\Sigma$ be the natural isomorphism, $\partial_c\in\{\Omega\Sigma A,A\}$ the connecting map associated
to the cone extension. The arguments of the proofs of \ref{theo:infsum_zero} and \ref{inverse}
show that, under the current hypothesis on $G$, $G(\Gamma A)=0$. Hence $\partial_c^G$ is an isomorphism.
On the other hand, the class of $\partial_c\tau$ in $kk(\Sigma A,A)$ is a natural isomorphism $\Omega\Sigma A\to A$.
Thus we must have
\[
\bar{G}(\alpha)=\partial_c^GG(\tau)\bar{G}(\Omega^n\alpha)G(\tau^{-1})(\partial_c^{G})^{-1}
\]
One checks that this rule does give a well-defined functor $kk\to \Ac$. To finish, we must show that $\bar{G}$ is homological. This amounts
to proving that $G$ maps \eqref{path:seq} to an exact sequence. By what we have already seen, the sequence
\[
\xymatrix{G(\Omega A)\ar[r]^{G(\Omega f)}&G(\Omega B)\ar[r]& G(P_f)\ar[r]& G(A)\ar[r]^{G(f)}&G(B)}
\]
is exact everywhere except perhaps at $G(A)$. Exactness at $G(A)$ follows by comparing this sequence with the path sequence
of the map $\Sigma \Omega f$.
\end{proof}
\sn
\begin{corollary} There is a natural action
\[
kk(A,B)\otimes G(A)\to G(B),\quad \alpha\otimes \xi\mapsto \alpha\cdot \xi=\bar{G}(\alpha)(\xi)\qed
\]
\end{corollary}

\begin{example} The functor $KH_0:\aha\to \rm{Ab}$ satisfies the hypothesis of Theorem \ref{univ:hom}. Hence we have a natural map
\begin{equation}\label{future_iso}
kk_0(\hal,A)\to \hom_{ab}(KH_0(\hal),KH_0(A))
\end{equation}
As $\hal$ is unital, there is a map $\ja :\Zz\to\hal$. Composing the map above with
\[
\hom_{ab}(KH_0(\hal),KH_0(A))\overset{\ja ^*}\longrightarrow\hom_{ab}(KH_0(\Zz),KH_0(A))\cong KH_0(A)
\]
we obtain a homomorphism $kk_0(\hal,A)\to KH_0(A)$. Applying this to $\Omega^nA$ we obtain
\begin{equation}\label{prebeta}
kk_*(\hal,A)\to KH_*(A)
\end{equation}
We will show in \ref{proof_main} that \eqref{prebeta} is an isomorphism.

\end{example}

\subsection{Comparison with Kassel's bivariant $K$-group}\label{subsec:compakas}

C. Kassel has constructed a bivariant abelian group $K(A,B)$, defined for any pair of central unital algebras $A,B\in\caha$ over a commutative
unital ground ring $\hal$, with unit-preserving structure maps. He defines $K(A,B)$ as the Grothendieck group of the exact category
$\rep_\hal(A,B)$ of those left $A\otimes_\hal B^{op}$-modules which are finitely generated projective as right $B$-modules. The exact
structure is given by the short exact sequences of bimodules
\begin{equation}\label{exact-kas}
0\to P'\to P\to P''\to 0.
\end{equation}
We point out that, as $P''$ is required to be projective, any exact sequence \eqref{exact-kas} is split in $B$-Mod,
and therefore in $H$-Mod, so that the added requirement of \cite[pp. 378]{kas} that
\eqref{exact-kas} remains exact upon tensoring with any $\hal$-module $V$ is automatic.

We wish to compare Kassel's group with $kk$, whenever both groups are defined.
A map $K(A,B)\to kk(A,B)$ is constructed as follows. If $P\in\rep(A,B)$, then there is an $\hal$-algebra
homomorphism $v:A\to\End_{B}(P)$. Choose a finitely generated projective $B$-module $Q$ such that $P\oplus Q=B^n$ and obtain
an inclusion $\End_{B}(P)\subset M_n(B)$. Composing with $v$, we obtain a homomorphism
\begin{equation}\label{up}
u(P):A\to M_n(B)
\end{equation}
which defines a class in $kk(A,B)$. Because $kk(A,?)$ is $M_2$-stable and split-exact, the construction above
preserves isomorphism classes. Moreover, a choice of a $B$-linear section of the map $P\to P''$ in \eqref{exact-kas} induces
a block matrix decomposition of $\End_B(P)$ under which
\[
u(P)(A)\subset\mathcal{U}:=\left[\begin{matrix}\End_B(P')&\hom_B(P'',P')\\0&\End_B(P'')\end{matrix}\right]\subset \End_B(P).
\]
The projection $\mathcal{U}\to \End_B(P')\oplus \End_B(P'')$ is a surjection with square-zero kernel, and therefore a $kk$-equivalence.
It follows that we have a group homomorphism
\begin{equation}\label{kas-to-kk}
u:K(A,B)\to kk(A,B).
\end{equation}
Kassel shows that the tensor product of bimodules defines a $\Zz$-bilinear composition product $K(A,B)\times K(B,C)\to K(A,C)$ via
$(P,Q)\mapsto P\otimes_B Q$. It is straightforward to check that \eqref{kas-to-kk} maps this composition to that of $kk$.
Furthermore, the tensor product over $\hal$ induces an external product (\cite[\S3]{kas})
\begin{equation}\label{kasext}
K(A_1,B_1)\otimes K(A_2,B_2)\to K(A_1\otimes_\hal A_2,B_1\otimes_\hal B_2).
\end{equation}
In general there is no analogue of this product in $kk$. However, if $\cU$ is the category of $\hal$-modules,
then we have the product \eqref{kkext}. It is straightfoward to show that \eqref{kas-to-kk} maps
\eqref{kasext} to \eqref{kkext} whenever both products are defined. Finally, Kassel obtains a Chern character (\cite[\S4]{kas})
\[
ch':K(A,B)\to HC^0(A,B).
\]
He defines $ch'$ as the composite of the trace map with the map \eqref{up} induces at the level of cyclic complexes. On the other
hand, we have seen in Example \ref{chern} that if $\hal$ is a field of characteristic zero we also have a map $ch:kk(A,B)\to HP^0(A,B)$.
It is immediate from the definitions that the canonical map $HC^0(A,B)\to HP^0(A,B)$ makes the following diagram commute
\[
\xymatrix{K(A,B)\ar[d]_{ch'}\ar[r]^u&kk(A,B)\ar[d]^{ch}\\
          HC^0(A,B)\ar[r] &HP^0(A,B).}
\]

\section{Some computations}\label{sec:compu}

In this section we provide some computations in the triangulated category $kk$. We first compute the $kk$-theory of an amalgamated product over a split
sub-algebra. Then we deal with graded algebras and nilpotent extensions. The last two subsections contain the fundamental theorem in $kk$-theory
and a computation for crossed products with $\Zz$.

\subsection{Amalgamated free products}\label{subsec:amalga}

Let $A,B$ be \hal-algebras containing a common subalgebra $C \subset A,B$. Furthermore, the existence
of retraction homomorphisms $\alpha: A \to C$ and $\beta: B \to C$ is assumed. In this situation, there is a natural map
\begin{equation}\label{amalgamap}
\vartheta: A \ast_C B \to A \oplus_C B,
\end{equation}
which is described by
\[a \mapsto (a,\alpha(a)) \quad \mbox{and} \quad b \mapsto (\beta(b),b)\]

This map was studied in \cite{cuntz_free}, where it was shown to
induce an isomorphism in topological algebra $KK$-theory. The same
argument, which we repeat for sake of completeness, applies in our
setting and shows that $j(\vartheta)$ is invertible.

\begin{theorem} \label{amalgamated} The map $\vartheta$ of \eqref{amalgamap} above is a $kk$-equivalence. \end{theorem}
\begin{proof}

First of all, we note that the map $\vartheta$ fits into a diagram of split extensions as follows.
\[ \xymatrix{D_1 \ar[d]^{\vartheta'} \ar[r] & A \ast_C B \ar[d]^{\vartheta}\ar[r]^{\alpha \ast \beta} & C \ar[d]^1 \\
D_2 \ar[r] & A \oplus_C B \ar[r] & C  }\]
Here, $D_2$ is equal to the subalgebra of those $(a,b)$, for which $\alpha(a)=\beta(b)=0$.
We are done if we can show that $j(\vartheta')$ is invertible. This is implied by the following Lemma. \end{proof}
\begin{lemma}
The map $\vartheta': D_1 \to D_2$ is a matrix-homotopy equivalence.
\end{lemma}
\begin{proof}
There is a homomorphism of \hal-algebras $\eta: D_2 \to M_2 D_1$ which maps
\[ (a,b) \mapsto  \left( \begin{array}{cc} a & 0 \\ 0 & b \end{array} \right) \in M_2 D_1\]
Note that the matrix $a \oplus b$ is indeed in $D_1$, since
$\alpha(a)=\beta(b)=0$, by our assumption that $(a,b) \in D_2$.
The composite $\eta  \vartheta': D_1 \to M_2(D_1)$ is induced by the map
$\phi: A \ast_C B \to M_2(A \ast_C B)$ which sends
\[a \mapsto a \oplus \alpha(a) \quad \mbox{and} \quad b \mapsto \beta(b) \oplus b.\]
Using a rotational homotopy, one shows that $\phi$ is homotopic to $1_{A \ast_C B} \oplus (\alpha \ast \beta)$.
Since $D_1$ is a direct summand of $A \ast_C B$ in $kk$,
on which $\alpha \ast \beta$ vanishes, we conclude that $j(\eta  \vartheta')$ is a $kk$-equivalence.

The other composition $M_2(\vartheta')  \eta: D_2 \to M_2(D_2)$ is described by
\[ (a,b) \mapsto (a,0) \oplus (0,b).\] It is clear, again by a rotational homotopy, that it is homotopic to the natural inclusion.
\end{proof}

\subsection{Graded algebras and nilpotent extensions}\label{subsec:gradednil}

\begin{theorem}\label{graded}
Let $A$ be a $\Nz$-graded \hal-algebra. The inclusion $A_0 \to A$ is a $kk$-equivalence.
\end{theorem}
\begin{proof} This follows from the well-known fact that the inclusion $A_0 \to A$
is an algebraic homotopy equivalence. The homotopy inverse is the projection $A \to A_0$.
Indeed, the map $A \to A[t]$ which sends an homogenous element
$a_n \in A_n$ to $a_n t^n$, is a homotopy between the composite $A \to A_0 \to A$ and the identity $1_A$.
\end{proof}

\begin{definition}
An $\hal$-algebra is {\rm nilpotent} if there exists a positive integer $m \in \Nz$ such that $A^m=0$, and $F$-nilpotent if it is nilpotent and, in addition,
all quotient maps $A^{2^n} \to A^{2^n}/A^{2^{n+1}}$ are
$F$-split surjections.
\end{definition}

\begin{remark}
In the particular case when $\cU={\rm Sets}$, and $F$ the forgetful functor, then every extension is $F$-split, and thus an
$F$-nilpotent algebra is the same as a nilpotent algebra.
\end{remark}

\begin{theorem}
Let $A$ be an $F$-nilpotent $\hal$-algebra. Then $A$ is $kk$-equivalent to zero.
\end{theorem}
\begin{proof} By induction the problem reduces to showing that square-zero algebras are $kk$-equivalent to zero. But any square-zero algebra $A$
is graded, with $A=A_1$, $A_n=0$ if $n\neq 1$. The theorem now follows from \eqref{graded}.\end{proof}

\subsection{The fundamental theorem}\label{subsec:fundatheo}

Let $A$ be an $\hal$-algebra. The fundamental theorem provides a computation of the Laurent polyomial ring $A[t^\pm]$ in the category $kk$.

\begin{theorem} \label{fundamental}
Let $A$ be an $\hal$-algebra. There is a $kk$-equivalence $A[t^\pm] \cong A \oplus \Sigma A$.
\end{theorem}
\begin{proof} First of all, it is clear that $A[t^\pm]$ splits into $A \oplus \sigma A$, by split-exactness.
We now analyze the following commutative diagram of $F$-split extensions
\[ \xymatrix{
M_\infty A \ar[r] \ar[d] & \tau_0 A \ar[r] \ar[d] & \sigma A \ar[d] \\
M_\infty A \ar[r]        & \Gamma A \ar[r] & \Sigma A. } \]

In order to show that the natural map $\sigma A \to \Sigma A$ is a $kk$-equivalence, it suffices to show that $\tau_0 A$ is
zero in $kk$. The result is therefore implied by the following lemma.
\end{proof}

\begin{lemma} \label{fundam_lemma}
Let $\Az$ be an abelian category, and $E:\aha\to \Az$  a split-exact, $M_\infty$-stable, homotopy invariant functor. Then $E(\tau_0) = 0$.
\end{lemma}
\begin{proof} Again, the proof is a more or less straightforward modification of the corresponding proof in the topological
setting, e.g. \cite[Prop. 8.2]{Cu_Weyl}. The strategy of the proof is to construct a quasi-homomorphism
from $\tau_0$ to $(M_\infty \tau)^{\Delta^1}$ which is the natural inclusion if evaluated at $t=0$ and which is zero at $t=1$. The result
then follows by Theorem \ref{quasi-homomorphisms}.
\mn
The quasi-homomorphism is constructed as follows. First of all, we define several maps from $\tau_0$ to $\tau \otimes \tau$, using the universal
properties of $\tau$.
\[\begin{array}{rclcrcl}
\psi(\alpha) &=& \beta \alpha^2 \otimes 1 &\quad& \psi(\beta)& =& \beta^2 \alpha \otimes 1 \\
\phi_1(\alpha) &=& \beta \alpha^2 \otimes 1 + e \otimes \alpha &\quad& \phi_1(\beta) &=& \beta^2 \alpha \otimes 1 + (1- \beta \alpha) \otimes \beta \\
\phi_2(\alpha) &=& \alpha \otimes 1  &\quad& \phi_2(\beta) &=& \beta \otimes 1 \\
\phi_3(\alpha) &=& \beta \alpha^2 \otimes 1 + e \otimes 1 &\quad& \phi_3(\beta) &=& \beta^2 \alpha \otimes 1 + (1-\beta \alpha) \otimes 1 \\
\end{array}\]
Furthermore, we define invertible elements $u_1$ and $u_2$ in $(\tau \otimes \tau)[t]$ as follows:
\[\begin{array}{lcl}
u_1 = \beta^2 \alpha^2 \otimes 1 + \left( \begin{array}{cc} 1 -t^2 \beta \alpha & (t^3-2t) \beta \\ t \alpha & 1-t^2 \end{array} \right) &\quad&
{u_1}^{-1} = \beta^2 \alpha^2 \otimes 1 + \left( \begin{array}{cc} 1 - t^2 \beta \alpha & t \beta \\ (t^3-2t) \alpha & 1-t^2 \end{array} \right) \\
u_2 = \beta^2 \alpha^2 \otimes 1 + \left( \begin{array}{cc} 1-t^2 & t^3-2t \\ t  & 1-t^2 \end{array} \right) &\quad&
{u_2}^{-1} = \beta^2 \alpha^2 \otimes 1 + \left( \begin{array}{cc} 1-t^2 & t \\ t^3-2t & 1-t^2 \end{array} \right). \end{array}\]
Here, we are implicitly using the inclusion $M_2 \otimes \tau[t] \to \tau \otimes \tau[t]$.
For $i=1,2$, we define homotopies $\Phi_i: \tau_0 \to \tau \otimes \tau[t]$ by the following formulas:
\[\begin{array}{rclcrcl}
\Phi_1(\alpha) &=& (\alpha \otimes 1)u_1  &\quad& \Phi_1(\beta) &=& {u_1}^{-1} (\beta \otimes 1) \\
\Phi_2(\alpha) &=& (\alpha \otimes 1)u_2  &\quad& \Phi_2(\beta) &=& {u_2}^{-1} (\beta \otimes 1).
\end{array}\]
One checks that $\ev_0  \Phi_1 = \ev_0  \Phi_2 =\phi_2$. On the other hand, $\ev_1  \Phi_1= \phi_1$ and
$\ev_1  \Phi_2 = \phi_3$.
Let $\Phi$ be the composite homotopy. It satisfies $\ev_0  \Phi = \phi_1$ and
$\ev_1  \Phi = \phi_3$.
Note that $\psi, \phi_1, \phi_2$ and $\phi_3$ agree modulo the ideal $M_\infty \otimes \tau$ and that $1- u_i$ is contained in $M_\infty \otimes \tau$.
All this together implies that there is a well-defined quasi-homomorphisms $(\Phi,\psi^{\Delta^1}): \tau_0 \to (M_\infty \otimes \tau)^{\Delta^1}$.
\mn
Now, obviously $\phi_1$ is the orthogonal sum of $\psi$ and the inclusion $\iota: \tau_0 \to M_\infty \otimes \tau$, hence
$E(\phi_1,\psi)= E(\iota)$ by Lemma \ref{split_sum}. On the other hand,
$\phi_3$ agrees with $\psi$ on $\tau_0$, so that $E(\phi_3,\psi)=0$. By homotopy invariance, it follows that the inclusion $\iota: \tau_0 \to M_\infty \otimes \tau$ induces the zero map.
By $M_\infty$-stability, this shows that
also the inclusion $\tau_0 \to \tau$ induces zero under $E$. This finishes the proof, using split-exactness once more.
\end{proof}

\begin{remark}
In the topological algebraic setting,  Theorem \ref{fundamental}
provides a proof of Bott periodicity, since in this case $\Omega$
and $\sigma$ coincide, due to the existence of a logarithm. In the
algebraic setting, this is no longer true.
\end{remark}

\subsection{Crossed products by $\Zz$}\label{subsec:cross}

Let $A\in\aha$ and $\sigma: A \to A$ be an automorphism of $\hal$-algebras.
We can form the crossed product algebra $A \rtimes_\sigma \Zz$. It is defined to be a twisted Laurent polynomial ring as follows.
As bimodule over $\hal$, $A \rtimes_\sigma \Zz$ is isomorphic to $A \otimes \Zz[t,t^{-1}]$. Multiplication is determined by the rule
\[
tat^{-1}=\sigma(a).
\]

The analysis of the $K$-theory of a crossed product in the operator algebraic setting was first carried out by
M. Pimsner and D. Voiculescu in \cite{PV}. They established
a $6$-term exact sequence, relating the operator $K$-theory of the crossed product $A \rtimes_\sigma \Zz$ to the $K$-theory of $A$. Hereafter, this sequence
was called Pimsner-Voiculescu exact sequence. Subsequently the argument of Pimsner and Voiculescu was extended to different
settings of topological algebras, e.g. Cuntz has extended the argument
to the setting of locally convex algebras in \cite{Cu_Weyl}.

The proof of the existence of a Pimsner-Voiculescu sequence in our setting can be taken verbatim from Section $14$ in \cite{Cu_Weyl}, subject to
suitable adjustments. We will recall the basic steps of the proof for convenience.

\begin{theorem}
Let $A\in\aha$ and $\sigma: A \to A$ an automorphism. There is an distinguished triangle
\[ \xymatrix{ \Omega A \ar[r] & A \ar[rr]^{1 -j(\sigma)} && A \ar[r]^(.4)\iota &  A \rtimes_\sigma \Zz.  } \]
\end{theorem}
\begin{proof}
Consider the subalgebra $\tau_\sigma$ of $\tau \otimes (A \rtimes_\sigma \Zz)$ which is generated by $1 \otimes A$, $\alpha \otimes t$ and $\beta \otimes t^{-1}$.
We get an extension
\[
0 \to M_\infty A \to  \tau_{\sigma} \to A \rtimes_\sigma \Zz \to 0.
\]
This $F$-split exact sequence induces the distinguished triangle above.
The proof consists of two steps. First of all, one shows that the natural map $\eta: A \to \tau_{\sigma}$ is an equivalence in $kk$-theory. Secondly, one
identifies the induced maps so that the conclusion of the theorem holds.

The map $\psi: \tau_\sigma \to \tau_\sigma$ which is defined by $a \mapsto (\beta \otimes 1) a (\alpha \otimes 1)$ is part of a
quasi-homomorphism $(1_{\tau_\sigma},\psi): \tau_\sigma \to \tau_\sigma \rhd M_\infty A$. It can be shown that $(1_{\sigma},\psi)$
is a two-sided inverse of $\eta: A \to \tau_\sigma$. Indeed, this is the content of \cite[Prop. 14.1]{Cu_Weyl}. The proof can be
taken from there, using the results from Lemma \ref{fundam_lemma}. The proof of \cite[Prop. 14.2]{Cu_Weyl} gives the identification
of the $kk$-elements as desired.
\end{proof}
\section{Comparison with $KH$}\label{sec:compa}

\subsection{Homotopy algebraic $K$-theory}\label{subsec:khdef}
Let $A$ be a unital ring. Put
\[
K(A)\deq BGl^+(A)
\]
for the connected $K$-theory space, and
\[
KV(A)\deq K(A^\Delta)
\]
for the realization of the simplicial space $[n]\mapsto
K(A^{\Delta^n})$. Set $K_n(A)=\pi_nK(A)$, $KV_n(A)=\pi_nKV(A)$
($n\ge 1$); these are respectively the Quillen and
Karoubi-Villamayor $K$-theory groups. Both theories preserve products. In
general if $G$ is any product preserving functor from unital rings
to abelian groups, and $A$ is any ring, we consider its
unitalization $\tilde{A}=A\oplus\Zz$ and set
\[
G(A)=\ker(G(\tilde{A})\to G(\Zz)).
\]
The Karoubi-Villamayor groups can be equivalently defined as
follows; for any ring $A$ we have
\[
KV_n(A)=\pi_n(BGl(A^\Delta))=\pi_{n-1}(Gl(A^\Delta)) \qquad (n\ge 1).
\]
In particular
\begin{align*}
KV_1(A)=&\pi_0(Gl(A^\Delta))\\
KV_{n+1}(A)=&KV_1(\Omega^nA)\qquad (n\ge 1).
\end{align*}

\sn
Next we consider the nonconnective $K$-theory spectrum $\Kt(A)$ of a unital ring $A$. To define it one uses the equivalence
of connected spaces
\begin{equation}\label{wagequiv}
K(A)\weq \Omega K(\Sigma A)_0
\end{equation}
and defines $\Kt(A)$ as the spectrum whose $n$th space is
\begin{equation}\label{ktn}
{}_n\Kt(A)\deq \Omega K(\Sigma^{n+1} A).
\end{equation}
The equivalence \eqref{wagequiv} induces a map
\begin{equation}\label{kvwagemap}
KV(A)\to \Omega KV(\Sigma A).
\end{equation}
Using this map one builds a nonconnective spectrum, the {\it homotopy $K$-theory} spectrum $\Kh(A)$, as follows. One
defines $\Kh(A)$ as the spectrum whose $n$th space is
\begin{align}\label{kv_kh_space}
{}_n\Kh(A)\deq &\Omega K(\Sigma^{n+1} A^{\Delta})\nonumber\\
=&\Omega K((\Sigma^{n+1}A)^{{\Delta}})\\
=&\Omega KV(\Sigma^{n+1}A).\nonumber
\end{align}
We put $K_n(A)=\pi_n(\Kt(A))$, $KH_n(A)=\pi_n\Kh(A)$ ($n\in \Zz$) for unital $A$, and extend these functors to the nounital case
as explained above. Note that, by definition of stable homotopy, if $A$ is unital we have
\begin{align}\label{formu_kh_kv}
KH_n(A)=&\colim_r\pi_{n+r}({}_r\Kh(A))\nonumber\\
       =&\colim_r KV_{n+r+1}(\Sigma^{r+1}A)\nonumber\\
       =&\colim_r KV_1(\Omega^{n+r}\Sigma^rA).
\end{align}
We remark that the formula $KH_n(A)=\colim_r KV_1(\Omega^{n+r}\Sigma^rA)$ is valid for not necessarily unital
$A$. This follows from the fact that $KV$ is split exact (\cite{KV1}).

\sn

In the next proposition we give an alternative formula for $KH_n(A)$. For this we need some natural
maps defined as follows. Let $A$ be any ring. Because $K_0(\Gamma A)=KV_1(\Gamma A)=0$, the surjection $K_1(\Sigma A)\to KV_1(\Sigma A)$
factors through $K_0(A)$, obtaining an epimorphism
\begin{equation}\label{k0kv1}
K_0(A)\fib KV_1(\Sigma A)
\end{equation}
On the other hand the loop extension gives a map $\partial:K_1(A)\to K_0(\Omega A)$ whose kernel is ${\rm Im} K_1(PA)$. A straightforward
verification shows that ${\rm Im} K_1(PA)=\ker (K_1(A)\to KV_1(A))$, so that $\partial(K_1(A))\cong KV_1(A)$. In particular
we have an injection
\begin{equation}\label{kv1k0}
KV_1(A)\hookrightarrow K_0(\Omega A)
\end{equation}
Applying \eqref{kv1k0} to $\Sigma A$ and composing with \eqref{k0kv1} we get a natural map $K_0(A)\to K_0(\Sigma \Omega A)$. Iterate this construction
and put
\[
KH'_n(A)\deq \colim_rK_0(\Sigma^r\Omega^{n+r}A)
\]
\begin{proposition}\label{kh=kh'}
Let $A$ be a ring. Then $KH_*(A)=KH'_*(A)$.
\end{proposition}
\begin{proof}
It suffices to show that the composite of \eqref{kv1k0} and of the natural transformation \eqref{k0kv1} applied to $\Omega A$ coincides
with the map used in the colimit \eqref{formu_kh_kv}. We have a commutative diagram
\[
\xymatrix{K_1(A)\ar[dr]^\iota\ar@{->>}[rr]&&KV_1(A)\ar[dddr]\ar[dr]^\iota&&K_2(\Sigma \tilde{A})\ar[dd]\ar[r]& KV_2(\Sigma\tilde{A})\ar[dd]\\
          &K_1(\tilde{A})\ar[ddrr]\ar[urrr]\ar[rr]&&KV_1(\tilde{A})\ar[urr]&\\
          &&&&K_1(\Sigma\widetilde{PA}:\Sigma\Omega A)\ar[r]&KV_1(\Sigma\Omega A)\\
          &&&K_0(\Omega A)\ar[ur]&&}
\]
Here $\iota$ is the inclusion $A\to \tilde{A}$. The horizontal maps are induced by the comparison map $K\to KV$,
and the oblique left to upper right maps are (induced by) the inverse of the edge map of the cone extension in $K$-theory.
The vertical maps in the back square as well as $K_1(\tilde{A})\to K_0(\Omega A)$ are induced by the extension
\[
\Omega A\to \widetilde{PA}\to \tilde{A}
\]
and $KV_1(A)\to K_0(\Omega A)$ is as defined above.
The homomorphism $KV_1(A)\to KV_1(\Omega\Sigma A)$ of \eqref{formu_kh_kv} is the composite shown in the diagram.
It is straightfoward to check that the composite map $K_0(\Omega A)\to KV_1(\Sigma \Omega A)$ of the
diagram is precisely the natural map \eqref{k0kv1} applied to $\Omega A$. The proposition follows by diagram chasing.
\end{proof}
The following corollary is due to Weibel. We give a new proof using \ref{kh=kh'}.
\begin{corollary}\label{easy}(\cite{weih})
$KH$ satisfies excision and is invariant under nilpotent extensions.
\end{corollary}
\begin{proof} Both assertions follow from the proposition above and properties of nonpositive $K$-theory.
\end{proof}

\subsection{Main theorem}\label{proof_main}

\begin{theorem} Let $\hal$ be a unital commutative ring. Assume $\aha=\caha$, and let $\cU$ be any of the underlying categories considered
in Section \ref{sec:conve}. Write $kk$ for the bivariant $K$-theory constructed from these data. Then the map \eqref{future_iso} induces an
isomorphism
\[
kk_*(\hal,A)\cong KH_*(A).
\]
\end{theorem}
\begin{proof}
Let $A$ be an $\hal$-algebra, $\ja :\Zz\to \hal$ the structure map.
Composition  with $\ja $ induces a bijection
\[
\hom_\aha(\hal,A)\weq\hom_\Ass(\Zz, A).
\]
This implies the following relation between the coproduct functors
$Q$ and $Q_\Zz$ computed respectively in $\aha$ and $\Ass$:
\[
Q(\hal)=\hal\otimes Q_\Zz(\Zz).
\]
It follows that $q(\hal)=\hal\otimes q_\Zz(\Zz)$. Write
$qq(\Zz,A)$ for the set of all quasi-homomorphisms
\[
(e_0,e_1):\Zz\rightrightarrows M_\infty \tilde{A}\vartriangleright
M_\infty A.
\]
Define a map
\begin{equation}\label{map_qq}
qq(\Zz,A)\to K_0(\tilde{A}), \qquad (e_0,e_1)\mapsto [e_0]-[e_1].
\end{equation}
The image of \eqref{map_qq} is precisely the subgroup
$K_0(A)=\ker(K_0(\tilde{A})\to K_0(\Zz))$. On the other hand, any
element of $qq(\Zz,A)$ gives rise to a ring homomorphism
$q_\Zz(\Zz)\to M_\infty A$, and thus to an $\hal$-algebra
homomorphism $q(\hal)\to M_\infty A$. One checks that the map
$qq(\Zz,A)\to kk_0(q(\hal),A)$ which sends a quasi-homomorphism to
the composite of the class of the corresponding map $q(\hal)\to A$
factors through a group homomorphism $K_0(A)\to kk_0(q(\hal),A)$.
Composing the latter with the inverse of $\delta:q(\hal)\to \hal$
we obtain
\begin{equation}\label{k0kk0}
\epsilon:K_0(A)\to kk_0(\hal,A).
\end{equation}
\sn \sn The map \eqref{k0kk0} induces
\[
\alpha:KH_0(A)=\colim_nK_0(\Sigma^n\Omega^nA)\longrightarrow
\colim_nkk_0(\hal,\Sigma^n\Omega^nA)=kk_0(\hal,A).
\]

Conversely, if $\partial_u$, $\partial_l$ and $\partial_c$ are the
connecting maps corresponding respectively to the universal, loop
(or any of its subdivided versions), and cone extensions, and
$e_0\in K_0(\Zz)$ the canonical generator, then the composite of
canonical maps
\[
\beta:kk_0(\hal,A)\to
\hom_{ab}(KH_0(\hal),KH_0(A))\overset{\ja ^*}\longrightarrow
\hom_{ab}(KH_0(\Zz),KH_0(A))\cong KH_0(A)
\]
sends the class of $\theta:J^n(\hal)\to M_rM_\infty A^{\sd^pS^n}$
to the image of $\ja (e_0)$ under
\[
\partial_l^{-n}\partial_c^nKH_0(\Sigma^n\theta)\partial_c^{-n}\partial_u^n:KH_0(\hal)\to KH_0(A).
\]
In particular, if $A$ is unital and $E\in M_\infty A$ an
idempotent, then the composite $\beta\alpha$ sends the class of
$E$ to its image in $KH_0(A)$. By construction of $\alpha$, this
is enough to prove that $\beta\alpha$ is the identity. To finish
the proof it suffices to show that $\alpha$ is onto. It is clear
that the class in $kk_0(\hal, A)$ of any map $\hal\to M_\infty A$
is in the image of $\alpha$. Let $n>0$ and $\theta:J^n(\hal)\to
M_rM_\infty B^{\sd^pS^n}$. Choose a map $\kappa:q(\hal)\to
M_\infty\Sigma^nJ^n(\hal)$ which represents the image of $\ja (e_0)$
under $\partial_c^{-n}\partial_u^{n}:K_0(\hal)\to
K_0(\Sigma^nJ^n(\hal))$. Note $\alpha(\ja (e_0))$ is the identity map
in $kk_0(\hal,\hal)$; hence we have the following equality in
$kk_0(\hal, \Sigma^nJ^n(\hal))$
\begin{equation}\label{aside}
j(\kappa)\cdot
j(\delta)^{-1}=\alpha((\partial_c^{-n}\partial_u^n)(\ja (e_0)))=\partial_c^{-n}\partial_u^n.
\end{equation}
In other words the composite of $\kappa$ with the inverse of
$\delta:q(\hal)\to \hal$ is the invertible excision element $\in
kk_0(\hal,\Sigma^nJ^n(\hal))$. Remark further that the following
diagram commutes in $kk$
\[
\xymatrix{\Sigma^n J^n(\hal)\ar[r]^{\Sigma^n\theta}&\Sigma^nB^{\sd^pS^n}\ar[d]^{\partial_l^{-n}\partial_c^n}\\
           \hal\ar[u]^{\partial_c^{-n}\partial_u^n}\ar[r]_\theta &B}.
\]
Put
$\theta'\deq\Sigma^n\theta\partial_c^{-n}\partial_u^{n}(\ja (e_0))$;
by \eqref{aside}, the composite of the upward arrow followed by
top row is $\epsilon(\theta')$. Let $\nu
:K_0(\Sigma^nB^{\sd^pS^n})\to KH_0(B)$ be composite of
$K_0(\Sigma^nB^{\sd^pS^n})\to KH_0(\Sigma^nB^{\sd^pS^n})$ followed
by $\partial_l^{-n}\partial_c^n$. Then
\[
\theta={\partial_l^{-n}\partial_c^n}\epsilon(\theta')=\alpha(\nu(\theta')).\qed
\]
\end{proof}
\begin{remark}\label{noncenkh}
Let $\hal$ be commutative, $A\in \gaha$. It follows from \eqref{noncentral} and the theorem above
that for the $kk$-groups $kk^b$
obtained from the choices $\aha=\gaha$ and $\cU=\hal\bimod$, we have
\[
kk_*^b(\hal,A)=KH_*(Z_\hal A).
\]
Here, as in \eqref{noncentral},
$Z_{\hal}A\subset A$ is the subring of those elements which commute with the action of $\hal$.
\end{remark}

\subsection{Bootstrap Category}\label{subsec:boot}

In this subsection we consider the natural question of whether the group $kk_*(A,B)$ can be computed in terms of $KH_*(A)$ and $KH_*(B)$.
For simplicity we restrict to the case $\hal=\Zz$. The construction of
a K\"unneth spectral sequence can in principle be carried out in the category $kk$, but the convergence of this spectral sequence is obstructed
by a simple observation. Although the symmetric ring spectrum $\kaha(\Zz)$ is a generator (in the triangulated sense, i.e. its suspensions
and de-suspensions form a generating set in the sense of \cite[ Definition $8.1.1$]{neeman}) of the triangulated category of $\kaha$-module
spectra, there is a priori no reason to expect that $\Zz$ is a generator for $kk$. In the similar case of $KK$-theory for $C^*$-algebras it is not known whether $\Cz$
is a generator for the corresponding triangulated category $KK$. Since the ingredients of the K\"unneth spectral sequence, i.e.\[E^{p,q}_2 = {\rm Ext}_{KH_*(\Zz)}^{p,q}(KH_*(A),KH_{*}(B)),\]
can only see $\kaha(\Zz)$-local information, the same is true for anything to which it could possibly converge.
It is fairly easy to see that the passage from $kk$ for the
triangulated category of $\kaha(\Zz)$-module spectra is precisely the
colocalization at the subcategory of $\Zz$-local objects, i.e. objects $A$ with $KH_*(A)=0$.

The preceding remarks explain that a suitable K\"unneth spectral sequence
should be constructed in the category of $\kaha(\Zz)$-module spectra and is in fact
already constructed in \cite{MMS} in a more general setting.
The remaining interesting and difficult question is the following:
Under what circumstances is
the map $kk(A,B) \to [\kaha(A),\kaha(B)]_{\kaha(\Zz)}$ an isomorphism?
The full subcategory formed by those $A$ for which the map above is an isomorphism for all $B$
is called the \textit{bootstrap category}. It is clear that the map will be an isomorphism
for all $B$, whenever $A$ is in the triangulated thick subcategory which is generated by $\Zz$.
Indeed, this follows just from an inductive construction of this subcategory and the $5$-lemma.
A further study of the size of the bootstrap category would require a better understanding
of the behaviour of $kk_*$ under limits and colimits in the category of $\aha$-algebras.
However, the authors have no results in that direction.

\end{document}